\documentclass[a4paper]{article}

\usepackage[utf8]{inputenc}
\usepackage[UKenglish]{babel}
\usepackage[T1]{fontenc}           	
\usepackage{amsmath}
\usepackage{amssymb}
\usepackage{amsfonts}          		
\usepackage[all]{xy}
\usepackage{amscd}
\usepackage{amsthm}			
\usepackage{graphics}
\usepackage{color}
\usepackage{mathtools}
\usepackage[a4paper]{geometry}
\usepackage{bbm}
\usepackage{mathrsfs}

\usepackage{lipsum}
\usepackage{indentfirst} 
\usepackage{latexsym} 
\usepackage{xcolor}
\usepackage{cite}
\usepackage[hidelinks]{hyperref}
\usepackage{algorithm}
\usepackage{algpseudocode}
\usepackage{enumitem}

\newtheorem{thm}{Theorem}[section]

\newtheorem{lem}[thm]{Lemma}
\newtheorem{prop}[thm]{Proposition}
\newtheorem{defn}{Definition}
\newtheorem{rem}{Remark}

\newtheorem{thm*}{Theorem}

\newtheorem{lem*}{Lemma}

\newcommand{\R}{\mathbb{R}}

\newcommand{\sgn}{\textnormal{sign}}
\newcommand{\abs}{\textnormal{abs}}

\newcommand{\pe}{\Pi_{\mathcal{E}}}
\newcommand{\one}{\mathbbm{1}}
\newcommand{\Pis}{\Pi_{\mathcal{S}}}
\newcommand{\range}{\textnormal{range}}
\newcommand{\io}{\mathcal{O}}

\DeclareMathOperator*{\argmin}{arg\,min}
\DeclareMathOperator*{\argmax}{arg\,max}

\numberwithin{figure}{section}
\numberwithin{table}{section}

\title{A low rank ODE for spectral clustering stability}
\author{\scshape{Nicola Guglielmi}\thanks{Gran Sasso Science Institute, L'Aquila, Italy (\href{mailto:nicola.guglielmi@gssi.it}{\tt nicola.guglielmi@gssi.it})}\scshape{ and Stefano Sicilia}\thanks{Gran Sasso Science Institute, L'Aquila, Italy (\href{mailto:stefano.sicilia@gssi.it}{\tt stefano.sicilia@gssi.it})}}
\date{}

\begin{document}
 \maketitle


 \begin{abstract}
  Spectral clustering is a well-known technique which identifies $k$ clusters in an undirected graph with weight matrix $W\in\mathbb{R}^{n\times n}$ by exploiting its graph Laplacian $L(W)$,
  whose eigenvalues $0=\lambda_1\leq \lambda_2 \leq \dots \leq \lambda_n$ and eigenvectors are related to the $k$ clusters. Since the  computation of $\lambda_{k+1}$ and $\lambda_k$ affects the reliability of this method, the $k$-th spectral gap $\lambda_{k+1}-\lambda_k$ is often considered as a stability indicator. This difference can be seen as an unstructured distance between $L(W)$ and an arbitrary symmetric matrix $L_\star$ with vanishing $k$-th spectral gap. A more appropriate structured distance to ambiguity such that $L_\star$ represents the Laplacian of a graph has been proposed by Andreotti et al. (2021). Slightly differently, we consider the objective functional
  $ F(\Delta)=\lambda_{k+1}\left(L(W+\Delta)\right)-\lambda_k\left(L(W+\Delta)\right)$,
  where $\Delta$ is a perturbation  such that $W+\Delta$ has non-negative entries and the same pattern of $W$. We look for an admissible perturbation $\Delta_\star$ of smallest Frobenius norm such that $F(\Delta_\star)=0$.  In order to solve this optimization problem, we exploit its low rank underlying structure. We formulate a rank-4 symmetric matrix ODE whose stationary points are the optimizers sought. The integration of this equation benefits from the low rank structure with a moderate computational effort and memory requirement, as it is shown in some illustrative numerical examples.
 \end{abstract}


 \section{Introduction}
 
 Clustering is the task of dividing a data set into $k$ communities such that members in the same groups are related. It is an unsupervised method in machine learning that discovers data groupings without the need of human intervention and its aim is to gain important insights from collected data. Spectral clustering (originating with Fiedler \cite{fiedler1973algebraic}) is a type of clustering that makes use of the Laplacian matrix of an undirected weighted graph to cluster its vertices into $k$ clusters. More precisely it performs a dimensionality reduction of the dataset and then it clusters in lower dimension. 
 
 The stability of this procedure is often associated with the spectral gap $g_k$, i.e. the difference between the $(k+1)$-st and $k$-th eigenvalues of the Laplacian. 
 When $g_k$ is not large, small perturbations may cause a coalescence of the two consecutive eigenvalues and may significantly change the clustering. Thus, according to the spectral gaps criteria, a suitable  number of clusters is the index of the largest spectral gap. This choice is also motivated by the fact that spectral gaps can be seen as an unstructured measure to ambiguity: up to a constant factor, $g_k$ represents the minimum of the Frobenius norm of the difference between the Laplacian and a symmetric matrix with coalescing $k$-th and $(k+1)$-st eigenvalues. Moreover the computation of the spectral gaps is not expensive.
 
 The problem of computing matrix stability distances arises in different fields of numerical linear algebra, where it is needed to compute verify the robustness of some data. Some examples are distance to singularity, matrix stability, measures in control theory,etc. (e.g \cite{guglielmi2015low,higham1988computing,higham2002computing,guglielmi2017nearest,kressner2015distance}).
 
 In this paper we introduce a structured measure to stability that takes into account the pattern of the weight matrix of the graph. In this way it is possible to achieve a result that is more appropriate than the one provided by the spectral gaps criteria. The distance considered here is similar to the one presented in \cite{andreotti2021measuring}, but in this case the different formulation allows to exploit the low rank underlying structure of the problem. This property leads to a significant memory savage thanks to a more efficient computation of the solution of the optimization problem given by the structured stability measure. 
 
 Our main objective is to describe in detail how to determine the new criteria for spectral clustering stability.
 We propose to compute the structured distance to ambiguity via a three-level approach, similar to the two-level approach of \cite{guglielmi2017matrix,guglielmi2015low}, which is divided in an \textit{inner iteration}, an \textit{outer iteration} and then a \textit{selection of $k$}, i.e. the best value for the number of clusters. The \textit{inner iteration} is the part of the algorithm that requires more effort: it consists in the solution of a non-convex eigenvalue optimization problem. In our method we see the optimizers of the problem as stationary points of a system of matrix ODEs whose size depends on the structural pattern of the weight matrix of the graph. Then, by generalizing the approach of \cite{guglielmi2022rank}, we define a rank-4 symmetric ODE whose stationary points are closely related to the full rank system and we integrate it until we reach convergence. When the $n\times n$ weight matrix has a number of nonzeros higher than $4n$, then it is more convenient to integrate the rank-4 ODE instead of the structured matrix ODE, with an important computational gain.
 
 The paper is organized as follows. In Section 2 we briefly describe the spectral clustering method and we illustrate how to measure its robustness by the introduction of a structured distance to ambiguity. In Section 3 we discuss how to solve the \textit{inner iteration} by means of a structured matrix ODE that is a gradient system. In Section 4 we exploit the low rank underlying structure of the gradient system to formulate a similar low rank ODE that is used to solve the \textit{inner iteration}. In Section 5 we describe the integration of the low rank ODE. Finally in Section 6 we present the numerical results of the algorithm in a few graphs with different features.
 
 \section{Distances to ambiguity for spectral clustering}
  
 Consider a graph $\mathscr{G}=(\mathscr{V},\mathscr{E},W)$, with $n$ vertices $\mathscr{V}$, edges $\mathscr{E}\subseteq \mathscr{V}\times \mathscr{V}$  and  weight matrix $W\in\R^{n\times n}$. Its Laplacian matrix is
 \[
  L=L(W)=\textnormal{diag}(W\mathbbm{1})-W, \qquad \mathbbm{1}=(1,\dots,1)^T.
 \]
 It is well-known that $L(W)$ is a symmetric and positive semi-definite matrix, so the spectral theorem ensures that its eigenvalues  
 \[
  \lambda_n(L(W))\geq \dots \geq \lambda_1(L(W))\geq 0
 \]
 are real non-negative and that their associated unit eigenvectors $x_n(L(W)),\dots,x_1(L(W))$ form an orthonormal basis of $\R^n$. Algorithm \ref{alg_spectral_clustering} shows how spectral clustering makes use of the spectrum of the Laplacian $L(W)$ (see \cite{von2007tutorial}) to partition the graph. The following result gives the theoretical reason behind the spectral clustering algorithm. 
 \begin{thm}
  \label{thm_kpart}
  Let $W\in \R^{n\times n}$ be the weight matrix of an undirected weighted graph $\mathscr{G}$ and denote by $L(W)$ its Laplacian. Then the number of the connected components $C_1,\dots,C_k$ of the graph equals the dimension of the kernel of $L(W)$. Moreover the eigenspace associated to the eigenvalue $0$ is spanned by the indicator vectors $\one_{C_1},\dots,\one_{C_k}$.
 \end{thm}

 \begin{algorithm}
  \caption{Unnormalized spectral clustering}
  \label{alg_spectral_clustering}
  \begin{description}
   \item[Input:] An undirected weighted graph $\mathscr{G}=(\mathscr{V},\mathscr{E},W)$ and the number of clusters $k$
   \item[Output:] Clusters $C_1,\dots,C_k$
  \end{description}
  \begin{algorithmic}[1]
   \State Find the $k$ smallest  eigenvalues $0=\lambda_1\leq \dots \leq \lambda_k$ of $L(W)$ and denote by $x_1,\dots,x_k\in\R^n$ the eigenvectors associated
   \State Build 
   \[
    X=\begin{pmatrix}
      & \vline & & \vline & & \vline & \\
      x_1 & \vline & x_2 & \vline & \cdots & \vline & x_k \\
      & \vline & & \vline & & \vline & \\
    \end{pmatrix}=
    \begin{pmatrix}
      & & r_1 & & \\
      \hline
      & & r_2 & & \\
      \hline
      & & \vdots & & \\
      \hline
      & & r_n & & \\
    \end{pmatrix}
    \]
   \State Associate $r_i$ with the $i$-th vertex of the graph 
   \item Cluster the points $r_1,\dots,r_n\in \R^k$ into $k$ clusters $C_1,\dots,C_k$
  \end{algorithmic}
 \end{algorithm}
   
 In order to evaluate the robustness of the clustering, it is crucial that the $k$ smallest eigenvalues of the Laplacian are not sensible to perturbations. Otherwise the eigenvectors associated may significantly change and hence the algorithm could lead to a completely different clustering. In this sense, spectral gaps provide a criteria to ensure a reasonable value of $k$. We can characterize them as the unstructured distance between the Laplacian and a symmetric matrix with vanishing spectral gap.
  
 \begin{thm}
  \label{thm_spectralgap}
  The $i$-th spectral gap $g_i=\lambda_{i+1}-\lambda_{i}$ is characterized as
  \[
   \frac{g_i}{\sqrt{2}}=\min\left\{ \|L(W)-\widehat{L}\|_F : \widehat{L}\in \textnormal{Sym}\left(\R^{n\times n}\right),\  \lambda_{i+1}(\widehat{L})=\lambda_i(\widehat{L})\right\},
  \]
  where $\textnormal{Sym}(\R^{n\times n})$ denotes the set of the symmetric real matrices.
  \begin{proof}
   See \cite[Theorem 3.1]{andreotti2021measuring}.
  \end{proof}
 \end{thm}
 
 However in the minimization problem of Theorem \ref{thm_spectralgap}, in general the optimizer is not the Laplacian of a graph, making this unstructured measure associated to the spectral gaps not so accurate. This motivates us to introduce a new stability measure that takes into account the structure of the weight matrix $W$, that is described by the sets
 \[
  \mathcal{S}=\left\{ A=(a_{ij})\in \R^{n\times n}: \  a_{ij}=0 \quad \forall (i,j)\notin \mathscr{E}\right\}, \qquad \mathcal{E}=\mathcal{S}\cap \textnormal{Sym}(\R^{n\times n}).
 \]
 We define the optimization problem 
 \begin{equation}
  \label{prob_dk}
  \Delta^{(k)}_\star=\argmin_{\Delta\in \mathcal{D}}\left\{ \|\Delta\|_F: \ \lambda_k(L(W+\Delta))=\lambda_{k+1}(L(W+\Delta)) \right\}, 
 \end{equation}
 where
 \[
  \mathcal{D}=\left\{ \Delta\in \mathcal{E}: \ W+\Delta\geq 0 \textnormal{ entrywise}\right\}
 \]
 is the set of all admissible perturbation that added to the weight matrix $W$ return a matrix with non-negative entries and with the same structure of $W$.
 The minimum of \eqref{prob_dk} 
 \[
  d_k(W)=\|\Delta^{(k)}_\star\|_F
 \]
 defines the $k$-th structured distance to ambiguity between $W$ and $W^{(k)}_\star:=W+\Delta^{(k)}_\star$. This new distance considered is similar to the one defined in \cite{andreotti2021measuring}, but it concerns a different geometry: in this framework we work with the Frobenius norm of the perturbation $\Delta$, instead of considering a unit normalization of $L(\Delta)$.  The reason behind this new choice is mostly practical, since in this way it is possible to exploit the low rank underlying properties of the problem by  the introduction of a rank-4 symmetric ODE.
 
 The approach presented relies on a three-level procedure:
 \begin{itemize}
  \item \textit{Inner iteration:} Given a perturbation size $\varepsilon>0$, we consider the non-negative objective functional
  \[
   F^{(k)}_\varepsilon(E)=\lambda_{k+1}\left(L(W+\varepsilon E)\right)-\lambda_k\left(L(W+\varepsilon E)\right),
  \]
  where the perturbation of $W$ is $\Delta=\varepsilon E$ with $\|E\|_F=1$. We look for a minimizer of the optimization problem
  \begin{equation}
   \label{prob_ii}
   E^{(k)}_\star(\varepsilon)=\argmin_{E\in \mathcal{D}_1} \ F^{(k)}_\varepsilon(E).
  \end{equation}
  \item \textit{Outer iteration:} We tune the parameter $\varepsilon$ to obtain the smallest value $\varepsilon_\star$ of the perturbation size such that the objective functional evaluated in the minimizer vanishes, that is
  \[
   F^{(k)}_{\varepsilon_\star}\left(E_\star^{(k)}(\varepsilon_\star)\right)=0.
  \]
  Then the optimizer of \eqref{prob_dk} would be $W^{(k)}_\star=W+\varepsilon_\star E_\star^{(k)}(\varepsilon_\star)$.
  \item \textit{Choice of $k$:} We repeat the procedure for all the values of $k\in [k_{\min},k_{\max}]$ and then select
 \[
  k_{\textnormal{opt}}(W)=\argmax_{k_{\min}\leq k \leq k_{\max}} d_{k}(W).
 \] 
 \end{itemize}
 
 \begin{rem}
  \label{not_semplif}
  Whenever the value of $k$ is fixed, we will omit it and we will denote, for brevity,
  \[
   x=x_{k+1}(L(W+\varepsilon E)), \quad y=x_k(L(W+\varepsilon E)), \quad \one=\sqrt{n}x_1(L(W+\varepsilon E)), 
  \]
  \[
   \lambda=\lambda_{k+1}(L(W+\varepsilon E)), \quad \mu=\lambda_k(L(W+\varepsilon E)), \quad 0=\lambda_1(L(W+\varepsilon E)).
  \]
 \end{rem}
 
 \begin{rem}
  \label{rem_sphere}
  Moreover, for any matrix set $\mathcal{A}$, we denote by  $\mathcal{A}_1$ its interjection with the unit Frobenius norm sphere
  \[
   \mathcal{A}_1=\left\{A\in \mathcal{A}:\|A\|_F=1\right\}.
  \]
 \end{rem}

 \section{A gradient system for the \textit{inner iteration}}
 
 In this section we describe an ODE based approach to solve the optimization problem \eqref{prob_ii} defined in the \textit{inner iteration}. We will consider as fixed parameters the value of $\varepsilon>0$ and a positive integer $k\in\{2,\dots,n-1\}$. We rewrite the perturbation as $\Delta=\varepsilon E$, with $\|E\|_F=1$ and we introduce a matrix path $E(t)\subseteq \mathcal{E}_1$ that represents the normalized perturbation of the weight matrix $W$.
 For the formulation of the ODE, we need the following time derivative formula for $F^{(k)}_\varepsilon(E(t)):=F_\varepsilon(E(t))$ (see  e.g. \cite{andreotti2019constrained} and \cite{andreotti2021measuring}).
 \begin{lem}
  \label{lem_derft}
  Let $E(t)$ be a differentiable path of matrices in $\mathcal{E}_1$ for $t\in [0,+\infty)$. Assume that, for a given $\varepsilon>0$, the eigenvalues $\lambda(t)=\lambda_{k+1}(L(W+\varepsilon E(t)))$ and $\mu(t)=\lambda_k (L(W+\varepsilon E(t)))$ are simple for all $t$. Then
  \[
   \frac{1}{\varepsilon}\frac{d}{dt}F_\varepsilon(E(t))=\langle G_\varepsilon(E(t)),\dot{E}(t)\rangle,
  \]
  where
  \[
   G=G_\varepsilon(E(t))=L^*\left(x(t)x(t)^T-y(t)y(t)^T\right)=\pe\left((x\bullet x- y\bullet y)\one^T-(xx^T-yy^T)\right)
  \]
  is the rescaled gradient of the objective functional $F_\varepsilon(E(t))$ and 
  \[
   \begin{matrix}
    L^*(W):& \R^{n \times n} & \rightarrow & \textnormal{Sym}(\R^{n\times n})\\ & W & \rightarrow & \pe\left(\textnormal{diagvec}(W)\one^T-W\right)
   \end{matrix}
  \] 
  is the adjoint of the Laplacian operator with respect to the inner Frobenius product $\langle X,Y\rangle=\textnormal{trace}(X^TY)$. Moreover $G_\varepsilon(E(t))\neq 0$ for all $t$.
 \end{lem}
 
 The gradient $G=G_\varepsilon(E(t))$ introduced in Lemma \ref{lem_derft} gives the steepest descent direction for minimizing the objective functional, without considering the constraint on the norm of $E$. The following result shows the best direction to follow in order to fulfill the unit norm condition, which can be rewritten as $\langle E,\dot{E}\rangle=0$.
 
 \begin{lem}
  \label{lem_opt_grad}
  Given $E\in \mathcal{E}_1$ and $G\in \mathcal{E}$, a solution of the optimization problem
  \begin{equation}
   \label{optE}
   \argmin_{Z\in \mathcal{E}_1, \  \langle Z,E \rangle=0} \langle G,Z \rangle
  \end{equation}
  is 
  \[
   \alpha Z_\star=-G+\langle G,E\rangle E,
  \]
  where $\alpha$ is the normalization parameter.
 \begin{proof}
  Let us consider the vectorized form of the matrices in $\R^{n^2}$. Then the Frobenius product in $\R^{n\times n}$ turns into the standard scalar product of $\R^{n^2}$ and the thesis is straightforward.
 \end{proof}

 \end{lem}
 
 Lemmas \ref{lem_derft} and \ref{lem_opt_grad} suggest to consider the matrix ordinary differential equation
 \begin{equation}
  \label{odeE}
  \dot{E}(t)=-G_\varepsilon(E(t))+\langle G_\varepsilon(E(t)),E(t) \rangle E(t),
 \end{equation}
 whose stationary points are zeros of the derivative of the objective functional $F_\varepsilon(E(t))$. Equation \eqref{odeE} is a gradient system for $F_\varepsilon(E(t))$, since along its trajectories
 \[
   \frac{d}{dt}F_\varepsilon(E(t))=\varepsilon\left(-\|G_\varepsilon(E(t))\|_F^2+(\langle G_\varepsilon(E(t)),E(t)\rangle)^2\right)\leq 0
 \]
 by means of the Cauchy-Schwartz inequality, which also implies that the derivative vanishes in $E_\star$ if and only if $E_\star$ is a stationary point of \eqref{odeE}. Thanks to the monotonicity property along the trajectories, an integration of this gradient system must lead to a stationary point $E_\star$. 
 
 The stationary point  $E_\star$ found belongs, by construction, to $\mathcal{E}_1$. Generally it also holds that $E_\star\in \mathcal{D}_1$ and in this case $E_\star$ is also a solution of \eqref{prob_ii}. However in the formulation of \eqref{odeE} it is not guaranteed that the stationary point found is an admissible perturbation of $W$ and hence a solution of optimization problem \eqref{prob_ii}. In order to ensure the admissibility of $E_\star$, we need to take into account the non-negative constraint $W+\varepsilon E_\star\geq 0$ componentwise.
 
 \subsection{Penalized gradient system}
 
 A possible way to impose that the path $E(t)$ is contained in $\mathcal{D}_1$ is by introducing the penalization term
 \[
  Q_\varepsilon(E)=\frac{1}{2}\left(\one^T(W+\varepsilon E)_-^2\one\right)=\frac{1}{2}\sum_{(i,j)\in \mathscr{E}}(w_{ij}+\varepsilon e_{ij})_-^2,
 \]
 where $(a)_-=\min(a,0)$ denotes the negative part of $a$. The new objective functional becomes
 \[
  F_{\varepsilon,c}(E)=F_\varepsilon(E)+c         Q_\varepsilon(E),
 \]
 where $c>0$ is the penalization size
 and the new optimization problem for the \textit{inner iteration} is
 \begin{equation}
  \label{prob_ii_pen}
  \argmin_{E\in \mathcal{E}_1} F_{\varepsilon,c}(E).
 \end{equation}
 In this way solutions of \eqref{prob_ii_pen} are forced to stay close to the set $\mathcal{D}$ if $c$ is big enough, in order to fulfill the non-negativity constraint of the weight matrix. 
 Now we show how the results for $F_\varepsilon(E)$ adapts to this new functional.
 \begin{lem}
  \label{lem_derft_pen}
  With the same hypothesis of Lemma \ref{lem_derft} it holds
  \[
   \frac{1}{\varepsilon} \frac{d}{dt}F_{\varepsilon,c}(E(t))=\langle G_{\varepsilon,c}(E(t)),\dot{E}(t)\rangle,
  \]
  where
  \[
   G_{\varepsilon,c}(E)=G_\varepsilon(E)+c(W+\varepsilon E)_-
  \]
  is the penalized gradient.
  \begin{proof}
   Since $E,\dot{E}\in\mathcal{E}$ are symmetric, we have
   \[
    \frac{d}{dt}Q_\varepsilon(E(t))=\sum_{(i,j)\in \mathscr{E}}\varepsilon \dot{e}_{ij}(t)(w_{ij}+\varepsilon e_{ij}(t))_-=\varepsilon\langle \dot{E}(t), (W+\varepsilon E(t))_-\rangle,
   \]
   where $\dot{E}(t)=(\dot{e}_{ij}(t))$.
   By repeating the same steps of the proof of Lemma \ref{lem_derft} we get the thesis.
  \end{proof}

 \end{lem}
 By replacing the gradient with the penalized gradient $G_{\varepsilon,c}(E)$, we obtain, as we did for equation \eqref{odeE}, the ODE
 \begin{equation}
  \label{odeE_pen}
  \dot{E}=-G_{\varepsilon,c}(E)+\langle G_{\varepsilon,c}(E),E \rangle E.
 \end{equation}
 In exactly the same way done for the non-penalized equation, we can show that equation \eqref{odeE_pen} is a gradient system whose stationary points are the only zeros of the derivative of $F_{\varepsilon,c}(E)$. Thus the trajectory $E(t)$ of equation \eqref{odeE_pen} is forced to stay close to $\mathcal{D}_1$, when $c$ is big enough, and hence the stationary points that will be reached are admissible solutions of problem \eqref{prob_dk} up to an error that is low if $c$ is huge.

 \section{A rank-4 symmetric equation}

 In this section we will consider a modified version of \eqref{prob_dk}, which does not take into account the non-negativity constraint of the set $\mathcal{D}$:
 \begin{equation}
  \label{prob_dk_nopen}
  \widetilde{W}^{(k)}_\star=\argmin_{\Delta\in \mathcal{E}}\left\{ \|\Delta\|_F: \ \lambda_k(L(W+\Delta))=\lambda_{k+1}(L(W+\Delta)) \right\},
 \end{equation}
 The introduction of this new problem is motivated by two different reasons. The former is that the presence of the non-negative constraint is difficult to insert in the low rank formulation that will be exposed. The latter is that in our experiments the violation of this constraint seems to be uncommon and hence generally $\widetilde{W}^{(k)}_\star$ and the solution of \eqref{prob_dk} coincide.
 
 Any solution of \eqref{prob_dk_nopen} that violates the constraint represents a weight matrix with one or more negative entries, which is not admissible; in these cases the low rank ODE approach is not suitable and we need to integrate the full rank system penalized \eqref{odeE_pen}.
 
 In case that the solution of \eqref{prob_dk} and \eqref{prob_dk_nopen} are the same, we propose a new matrix ODE whose aim is to exploit the underlying low rank property of the problem and to solve more efficiently the \textit{inner iteration}. This derivation is a generalization of the rank-1 ODE based approach exhibited in \cite{guglielmi2022rank}.
 
 \subsection{Formulation of the low rank symmetric ODE}
 
 We introduce two low rank matrices $N$ and $R$, whose formulation depends on the matrix $E$, on the perturbation size $\varepsilon$ and on the fixed positive integer $k$ (which will be omitted):
 \[
  N=N_\varepsilon(E)=z \one^T-xx^T+yy^T, \qquad 
  R=R_\varepsilon(E)=\frac{N+N^T}{2}=\frac{z\one^T+\one z^T}{2}-xx^T+yy^T,
 \]
 where
 \[
  z=x\bullet x- y\bullet y
 \]
 is the vector of entries $z_i=x_i^2-y_i^2$ and $\bullet$ denotes the componentwise product. 
 We observe that the gradient can be rewritten as
 \[
  G=G_\varepsilon(E)=\pe(N_\varepsilon(E))=\Pis(R_\varepsilon(E)),
 \]
 which means that $G$ is the projection onto the pattern given by $\mathcal{S}$ of the low rank symmetric matrix $R$.
 
 \begin{rem}
  \label{rem_z}
  Since $x$ and $y$ have unit norm, we observe that 
  \[
   z^T \one=\sum_{i=1}^n (x_i^2-y_i^2)=1-1=0,
  \]
  which means that $\one$ and $z$ are orthogonal.  The vectors $x,y$ and $z$ are generally linear independent, but it may happen that they are not; however this seems to be a very exceptional condition.
  In the following we will assume that $x,y$ and $z$ are linearly independent, which implies, that the matrix $N$ has rank $3$ and hence $R$ has rank $4$.
 \end{rem} 
 Throughout this section and in the next ones we will use a particular type of decomposition of a low rank matrix symmetric matrix, which mixes the properties of the SVD and the spectral decompositions.
 \begin{defn}
  Let $Y\in \textnormal{Sym}(\R^{n\times n})\cap \mathcal{M}_r$ be a symmetric rank-$r$ matrix, where $\mathcal{M}_r$ denotes the rank-$r$ manifold. Then a singular values symmetric decomposition (SVSD) is
  \[
    Y=USU^T,
  \]
  where $U\in \R^{n\times r}$ has full rank and orthonormal columns and $S\in \textnormal{Sym}(\R^{r\times r})$ is invertible.
 \end{defn}

 \begin{rem}
  \label{rem_R}
  The matrix $R$ can be rewritten in the form
  \[
  R=\left(\frac{z+\one}{2}\right)\left(\frac{z+\one}{2}\right)^T-\left(\frac{z-\one}{2}\right)\left(\frac{z-\one}{2}\right)^T-xx^T+yy^T=R_U R_S R_U^T,
 \]
 where
 \[
  R_U=\begin{pmatrix}
   z+\one & z-\one & x & y
  \end{pmatrix}, \qquad 
  R_S=\begin{pmatrix}
   \frac{1}{4} & 0 & 0 & 0 \\
   0 & -\frac{1}{4} & 0 & 0 \\
   0 & 0 & -1 & 0 \\
   0 & 0 & 0 & 1 \\
  \end{pmatrix}.
 \]
 Thus, by means of a QR decomposition of $R_U$ it is possible to obtain an SVSD decomposition 
 \[
  R=B\Lambda B^T
 \]
 where the orthogonal $n$-by-$4$ matrix $B$ depends smoothly on $\one,x,y$ and $z$ and $\Lambda$ is a $4$-by-$4$ invertible symmetric matrix.
 \end{rem}
  
 Solutions of \eqref{odeE}
 can be rewritten as $E=\Pis Z$, where $Z$ solves
 \begin{equation}
  \label{odeZ}
  \dot{Z}=-R_\varepsilon(E)+\langle R_\varepsilon(E),E \rangle Z,
 \end{equation}
 and we recall $G_\varepsilon(E)=\Pis R_\varepsilon(E)$. We take inspiration from equation \eqref{odeZ} and we consider the ODE in the rank-4 manifold $\mathcal{M}_4$
 \begin{equation}
  \label{odeY}
  \dot{Y}=-P_Y R_\varepsilon(E)+\eta Y, \qquad \eta=\langle P_Y R_\varepsilon(E),E \rangle, \qquad E=\Pis Y,
 \end{equation}
 where $P_Y$ is the orthogonal projection, with respect to the Frobenius inner product, onto the tangent space $\mathcal{T}_Y \mathcal{M}_4$ at $Y$. If $Y=V_1\Sigma V_2^T$ is an SVD decomposition of $Y$, then the expression of $P_Y$ is given by the formula (see \cite{koch2007dynamical})
 \[
  P_Y A=A-(I_n-V_1V_1^T)A(I_n-V_2V_2^T),
 \]
 where $I_n$ denotes the $n$-by-$n$ identity matrix.
 
 \begin{rem}
  \label{svdrem}
  If $Y=USU^T$ is an SVSD decomposition of a rank-$r$ symmetric matrix $Y$, then the projection $P_Y$ 
  onto the tangent space $\mathcal{T}_Y\mathcal{M}_r$ can be rewritten as
  \[
   P_Y A=A-(I_n-UU^T)A(I_n-UU^T), \qquad \forall A\in \R^{n\times n}.
  \]
  To prove this fact we can consider a spectral decomposition $S=QDQ^T$, where $Q$ is orthogonal and $D$ is invertible and diagonal with elements ordered increasingly in absolute value. Then the associated singular values decomposition is
  \[
   Y=UQ\abs(D)\sgn(D)Q^TU^T, \qquad V_1=UQ, \quad \Sigma=\abs(D), \quad V_2=UQ\sgn(D)
  \]
  where $\abs(D)$ and $\sgn(D)$ are the matrices with absolute values and sign (respectively) of the diagonal elements and it holds $V_1V_1^T=UU^T=V_2V_2^T$.
 \end{rem}
 
 The following proposition shows the main structural features of the solution $Y(t)$ of \eqref{odeY}.
 
 \begin{prop}
  \label{prop_Y}
  Let $Y(t)$ be a solution of equation \eqref{odeY} for $t\in[0,+\infty)$ with starting value $Y(0)=Y_0\in \textnormal{Sym}(\R^{n\times n})\cap \mathcal{M}_4$. Then $Y(t)\in\textnormal{Sym}(\R^{n\times n})\cap \mathcal{M}_4$ for all $t$.
  
  Moreover, if $\|\Pis Y_0\|_F=1$, then $\|\Pis Y(t)\|_F=1$ for all $t$.
  \begin{proof}
   We show the properties by means of a differential argument that holds for $t$ close to 0 and that then extends to all $t$. Since $P_Y Y=Y$ and $P_Y R\in \textnormal{Sym}(\R^{n\times n})$, the derivative of the solution $Y(t)$ of \eqref{odeY} is
   \[
    \dot{Y}=-P_Y R+\langle P_Y R,E \rangle Y=P_Y\left(-R+\langle P_Y R,E \rangle Y\right)\in \mathcal{T}_Y \mathcal{M}_4\cap \textnormal{Sym}(\R^{n\times n}),
   \]
   which means that, for all $t$, the matrix $Y(t)$ lays in the rank-4 manifold and is symmetric.
   Similarly, if we define the matrix
   \[
    E(t)=\Pis Y(t)
   \]
   and we assume that $\|E(0)\|_F=\|\Pis Y_0\|_F=1$, then its unitary Frobenius norm is preserved, 
   \[
    \frac{d}{dt}\|E(t)\|_F^2=\langle \dot{E}(t),E(t)\rangle=-\langle P_Y R_\varepsilon(E),E\rangle+\eta \langle Y,E\rangle=-\langle P_Y R_\varepsilon(E),E\rangle+\eta\|E\|_F^2=0.
   \]
  \end{proof}
 \end{prop}
 
 In the next paragraphs we investigate how the solution of \eqref{odeY} is related to that of \eqref{odeE}. More precisely we are interested in their stationary points and in the monotonicity property of the low rank system, which are crucial for the implementation of the \textit{inner iteration}. If these properties are shared between the equations, then it would be possible to integrate the low rank ODE instead of the original ODE.

 \subsection{Comparison of the stationary points}
 
 It turns out that equations \eqref{odeE} and \eqref{odeY}, under non-degeneracy conditions, share the same stationary points.
 Before stating this result, we need the following technical lemma.
 
 \begin{lem}
  \label{lem_R_RUU}
  Let $R,Y\in \R^{n\times n}$  be two symmetric matrices with the same range of dimension $r$. Consider an SVSD decomposition $Y=USU^T$. Then
  \[
   R=RUU^T.
  \]
  \begin{proof}
   It is well known that the projector from $\R^n$ onto the range of $U$ with respect to the standard scalar product is $\Pi_U=U(U^TU)^{-1}U^T=UU^T$. Indeed for all $v\in \R^n$ and for all $x\in \R^r$
   \[
    \Pi_U(v)\in \range(U), \qquad \langle U(U^TU)^{-1}U^Tv,Ux\rangle=v^TU(U^TU)^{-1}U^TUx=\langle v,Ux\rangle.
   \]
   Since $R$ and $Y$ share the same range and they are symmetric, it holds
   \[
    \range(Y)=\range(R) \Rightarrow \range(Y)^\perp=\range(R)^\perp \Rightarrow \ker(Y^T)=\ker(R^T) \Rightarrow \ker(Y)=\ker(R).
   \]
   Moreover also $U$ and $Y$ share the same range: indeed these spaces have both dimension $r$ and it is straightforward that the range of $Y$ is included in the range of $U$.
   Hence any vector $v\in\R^n=\ker(Y)\oplus \range(Y)$, can be decomposed as
   \begin{equation}
    \label{lemdec}
    v=UU^Tv+y,
   \end{equation}
   for some $y\in \ker(Y)=\ker(R)$. 
   Finally applying $R$ to both sides of $\eqref{lemdec}$ yields the thesis
   \[
    Rv=RUU^Tv, \qquad \forall v\in \R^n.
   \]
  \end{proof}
 \end{lem}
 
 Now we can prove the main result on the stationary points of equations \eqref{odeE} and \eqref{odeY}. 
 
 \begin{thm}
  \label{thm_stat}
  Consider the two matrix ordinary differential equations
  \begin{equation}
   \label{Eode}
   \dot{E}=-G_\varepsilon(E)+\langle G_\varepsilon(E),E \rangle E,
  \end{equation}
  \begin{equation}
   \label{Yode}
   \dot{Y}=-P_Y R_\varepsilon(E)+\langle P_Y R_\varepsilon(E),E \rangle Y,
  \end{equation}
  \begin{enumerate}
   \item Let $E_\star\in \mathcal{E}_1$ of unit Frobenius norm be a stationary point of \eqref{Eode}. Then $E_\star=\Pis Y_\star$ for a certain symmetric matrix $Y_\star\in \mathcal{M}_4$ that is a stationary point of \eqref{Yode}.
   \item Conversely, let $Y_\star\in \mathcal{M}_4$ be a symmetric  stationary point of \eqref{Yode} such that $E_\star=\Pis Y_\star$ has unit Frobenius norm and $P_{Y_\star} R_\star\neq 0$, where $R_\star=R_\varepsilon(E_\star)$. Then $P_{Y_\star} R_\star=R_\star$, $Y_\star$ is a nonzero real multiple of $R_\star$ and $E_\star$ is a stationary point of \eqref{Eode}.
  \end{enumerate}
  \begin{proof}
   \begin{enumerate}
    \item Let $E_\star$ of unit Frobenius norm be a stationary point of \eqref{Eode}. Then, by Lemma \ref{lem_derft}, $G_\varepsilon(E_\star)\neq 0$ and there exists $\nu \in \R\setminus\{0\}$ such that 
    \[
     E_\star=\nu^{-1} G_\star=\nu^{-1} \Pis R_\star\neq 0.
    \]
    Defining the symmetric matrix  $Y_\star:=\nu^{-1} R_\star \in \mathcal{M}_4$ yields $E_\star=\Pis Y_\star$ and
    \[
     P_{Y_\star} R_\star=\nu P_{Y_\star} Y_\star=\nu Y_\star=R_\star.
    \]
    Moreover
    \[
     \langle P_{Y_\star} R_\star,E_\star \rangle=\langle R_\star,E_\star \rangle=\langle \Pis R_\star,E_\star \rangle=\nu \|E_\star\|_F^2=\nu.
    \]
    Thus the left hand side of \eqref{Yode} is
    \[
     -P_{Y_\star} R_\star+\langle P_{Y_\star} R_\star,E_\star \rangle Y_\star=-R_\star+\nu Y_\star=0,
    \]
    which means that $Y_\star$ is a stationary point of \eqref{Yode}.
    \item We begin by showing that $Y_\star$ is a nonzero real multiple of $R_\star$.
    The hypothesis yields $Y_\star=\nu^{-1} P_{Y_\star} R_\star\neq 0$ for some $\nu\in \R\setminus\{0\}$, that is
    \begin{equation}
     \label{RW}
     R_\star=\nu Y_\star+W,
    \end{equation}
    where $W\in \R^{n\times n}$ satisfies $P_{Y_\star} W=0$.
    Let $Y_\star=U_\star S_\star U_\star ^T$ 
    be an SVSD decomposition of $Y_\star$. Then
    \[
     W=(I-U_\star U_\star^T)W(I-U_\star U_\star^T)
    \]
    and equation \eqref{RW} becomes
    \[
     R_\star=\nu U_\star S_\star U_\star^T+(I-U_\star U_\star^T)W(I-U_\star U_\star^T).
    \]
    By multiplying from the right by $U_\star$ we get
    \[
     R_\star U_\star=\nu U_\star S_\star,
    \]
    which means, as shown in lemma \ref{lem_R_RUU}, that $R_\star,Y_\star$ and $U_\star$ have the same range. Then  
    \[
     R_\star=R_\star U_\star U_\star^T=\nu U_\star S_\star U_\star^T U_\star U_\star^T+(I-U_\star U_\star^T)W(I-U_\star U_\star^T)U_\star U_\star^T=\nu Y_\star,
    \]
    which shows that $Y_\star$ is a nonzero multiple of $R_\star$. Moreover $E_\star=\Pis Y_\star$ is a stationary point of \eqref{Eode} since
    \[
     \langle P_{Y_\star} R_\star,E_\star \rangle=\langle R_\star, E_\star \rangle=\nu\langle Y_\star,E_\star \rangle = \nu\|E_\star\|_F^2=\nu
    \]
    and hence
    \[
     -G_\varepsilon(E_\star)+\langle G_\varepsilon(E_\star),E_\star \rangle E_\star=-\Pis R_\star+\langle R_\star,E_\star\rangle \Pis Y_\star=-\nu \Pis Y_\star+\nu \Pis Y_\star=0.
    \]
   \end{enumerate}
  \end{proof}
  
 \end{thm}

 \subsection{Local convergence to the stationary points of the rank-4 ODE}
 
 Theorem \ref{thm_stat} ensures that the original and the low rank ODEs share the same stationary points. Now we are interested in understanding whether the integration of \eqref{odeY} leads to at least one of the local minima (i.e. the stationary points of the low rank ODE) or not. This convergence is always guaranteed for equation \eqref{odeE}, since it is a gradient system, but unfortunately equation \eqref{odeY} is not a gradient system and the monotonicity property of the functional may not hold.
 
 However, provided a suitable starting value for the integration of the ODE sufficiently close to a local minimum, the low rank ODE turns out to be close to a gradient system. The following key lemma shows the main reason behind this fact. 
 
 \begin{lem}
  \label{lem_key}
  Let $Y_\star\in \mathcal{M}_4\cap \textnormal{Sym}(\R^{n\times n})$ be a stationary point of the rank-$4$ ODE \eqref{odeY} such that $E_\star=\Pis Y_\star\in \mathcal{E}_1$ and $P_{Y_\star}R(E_\star)\neq 0$. Then there exists $\delta_\star>0$ such that for all $\hat{Y}\in \mathcal{M}_4$ that satisfies
  \[
   \|\hat{Y}-Y_\star\|_F=\delta<\delta_\star, \qquad \hat{E}=\Pis \hat{Y}\in \mathcal{E}_1,
  \]
  it holds
  \[
   \|P_{\hat{Y}} R_\varepsilon(\hat{E})-R_\varepsilon(\hat{E})\|_F\leq C\delta^2,
  \]
  where $C$ is a positive constant independent of $\delta$. Thus, $R_\varepsilon(\hat{E}(t))$ and $P_{\hat{Y}} R_\varepsilon(\hat{E}(t))$ coincide, up to quadratic terms, in a right-neighborhood of $0$.
  \begin{proof}
   Let $Y_\star$ be a stationary point of \eqref{odeY} that fulfills the hypothesis. Then, as shown in theorem \ref{thm_stat}, it holds
   \[
    Y_\star=U_\star S_\star U_\star^T=\nu^{-1} R_\star=B_\star \Lambda_\star B_\star^T,
   \]
   where $U_\star=B_\star\in \R^{n\times 4}$ and $S_\star=\nu^{-1} \Lambda_\star\in \R^{4\times 4}$.
   We introduce the matrix paths
   \[
    \widetilde{Y}(\tau)=\widetilde{U}(\tau) \widetilde{S}(\tau) \widetilde{U}(\tau)^T, \qquad \widetilde{R}(\tau)=\widetilde{B}(\tau)\widetilde{\Lambda}(\tau) \widetilde{B}(\tau)^T, \qquad \tau\in[0,\delta]
   \]
   which are guaranteed to be smooth and such that 
   \[
    Y_\star=\widetilde{Y}(0), \qquad R_\star=\widetilde{R}(0), \qquad \hat{Y}=\widetilde{Y}(\delta), \qquad R_\varepsilon(\hat{E})=\widetilde{R}(\delta).
   \]
   Since any matrix $\hat{Y}$ that satisfies the hypothesis can be written as $\hat{Y}=\widetilde{Y}(\delta)$, for instance
   \[
    \widetilde{Y}(\tau)=\hat{Y}+\frac{\tau-\delta}{\delta}\left(\hat{Y}-Y_\star\right),
   \]
   it is enough to study these paths in order to conclude.
   
   We will denote, for brevity, by $U,B,\Lambda$ and later $\dot{U},\dot{B},\dot{\Lambda}$ the associated function (equipped with the $\sim$) evaluated at $\tau=0$.
   The derivatives of $\widetilde{B}$ and of the other matrix functions, are well defined in a right-neighborhood of $\tau=0$. Indeed the explicit formulas (see \cite{meyer1988derivatives} for more details) for the eigenvectors' derivatives  
   \[
    \frac{d}{d\tau}\tilde{x}(\tau)=-\varepsilon \left(L\left(W+\varepsilon E\left(\widetilde{Y}(\tau)\right)\right)-\lambda I\right)^\sharp L\left(\dot{E}\left(\widetilde{Y}(\tau)\right)\right)x(\tau)
   \]
   \[
    \frac{d}{d\tau}\tilde{y}(\tau)=-\varepsilon \left(L\left(W+\varepsilon E\left(\widetilde{Y}(\tau)\right)\right)-\mu I\right)^\sharp L\left(\dot{E}\left(\widetilde{Y}(\tau)\right)\right)y(\tau)
   \]
   show that the first derivative of $\tilde{x}(\tau)$ and $\tilde{y}(\tau)$ are bounded by exploiting the group inverse (here denoted by $\sharp$), which in this case coincides with the more familiar Moore-Penrose pseudo-inverse. Hence, by means of remark \ref{rem_R}, also $\dot{U},\dot{B}$ and $\dot{\Lambda}$ are well defined and this allows to expand until the first order the matrices $\widetilde{R}$ and $P_{\tilde{Y}} \widetilde{R}$ for $0\leq\tau\leq \delta$. Recalling that $U(0)=U_\star=B(0)=B_\star$ yields
   \[
    P_{\tilde{Y}(\tau)} \widetilde{R}(\tau)=\widetilde{R}(\tau)\widetilde{U}(\tau)\widetilde{U}(\tau)^T+\widetilde{U}(\tau)\widetilde{U}(\tau)^T\widetilde{R}(\tau)-\widetilde{U}(\tau)\widetilde{U}(\tau)^T\widetilde{R}(\tau)\widetilde{U}(\tau)\widetilde{U}(\tau)^T=
   \]
   \[
    =B\Lambda B^T+\tau \left(\dot{B}\Lambda B^TUU^T+B\dot{\Lambda}B^TUU^T+B\Lambda \dot{B}^TUU^T+B\Lambda B^T\dot{U}U^T+B\Lambda B^TU\dot{U}^T\right)+
   \]
   \[
    +\tau\left(\dot{U}U^TB\Lambda B^T+U\dot{U}^TB\Lambda B^T+UU^T\dot{B}\Lambda B^T+UU^TB\dot{\Lambda}B^T+UU^TB\Lambda \dot{B}^T\right)+
   \]
   \[
    +\tau\left(-\dot{U}U^TB\Lambda B^TUU^T-U\dot{U}^TB\Lambda B^TUU^T-UU^T\dot{B}\Lambda B^TUU^T-UU^TB\dot{\Lambda} B^TUU^T\right)+
   \]
   \[
    +\tau\left(-UU^TB\Lambda \dot{B}^TUU^T-UU^TB\Lambda B^T\dot{U}U^T-UU^TB\Lambda B^TU\dot{U}^T\right)+\io(\tau^2)=
   \]
   \[
    =B\Lambda B^T+\tau\left(\dot{B}\Lambda U^T+U\dot{\Lambda}U^T+U\Lambda \dot{B}^TUU^T+U\Lambda U^T \dot{U}U^T+U\Lambda \dot{U}^T\right)+
   \]
   \[
    +\tau\left(\dot{U}\Lambda U^T+U\dot{U}^TU\Lambda U^T+UU^T\dot{B}\Lambda U^T+U\dot{\Lambda}U^T+U\Lambda \dot{B}^T\right)+\io(\tau^2)+
   \]
   \[
    +\tau\left(-\dot{U}\Lambda U^T-U\dot{U}^TU\Lambda U^T-UU^T\dot{B}\Lambda U^T-U\dot{\Lambda}U^T-U\Lambda \dot{B}^TUU^T-U\Lambda U^T\dot{U}U^T-U\Lambda \dot{U}^T\right)=
   \]
   \[
    =B\Lambda B^T+\tau\left(\dot{B}\Lambda U^T+U\dot{\Lambda}U^T+U\Lambda \dot{B}^T\right)+\io(\tau^2),
   \]
   while
   \[
    \widetilde{R}(\tau)=B\Lambda B^T+\tau(\dot{B}\Lambda B^T+B\dot{\Lambda} B^T+B\Lambda \dot{B}^T)+\io(\tau^2)
   \]
   and the thesis is straightforward.
   
  \end{proof}
  
 \end{lem}
 
 Now we are ready to state and prove the local convergence result to a strong local minimum.
 
 \begin{thm}
  \label{thm_loc_conv}
  Let $Y_\star \in \mathcal{M}_4\cap \textnormal{Sym}(\R^{n\times n})$ be a stationary point of the projected differential equation \eqref{odeY} such that $E_\star=\Pis Y_\star\in \mathcal{S}_1$ and $P_{Y_\star}R_\varepsilon(E_\star)\neq 0$. Suppose that $E_\star$ is a strong local minimum of the functional $F_\varepsilon$ on $\mathcal{S}_1$ and assume that 
  \[
   \Pis \Big \lvert_{\mathcal{M}_4}:\mathcal{M}_4\rightarrow \Pis(\mathcal{M}_4)\subseteq \mathcal{S}
  \]
  is a diffeomorphism. Then, for an initial datum $Y(0)$ sufficiently close to $Y_\star$, the solution $Y(t)$ of \eqref{odeY} converges to $Y_\star$ exponentially as $t\rightarrow +\infty$. Moreover $F_\varepsilon(\Pis Y(t))$ decreases monotonically with $t$ and converges exponentially to the local minimum value $F(E_\star)$ as $t\rightarrow +\infty$.
  \begin{proof}
   Thanks to the assumption that $\Pis$ is a diffeomorphism between $\mathcal{M}_4$ and its image, the differential equation \eqref{odeY} is equivalent to
   \[
    \dot{E}=-\Pis P_Y R_\varepsilon(E)+\langle \Pis P_Y G_\varepsilon(E),E \rangle E
   \]
   and hence Lemma \ref{lem_key} implies
   \[
    \dot{E}=-\Pis R_\varepsilon(E)+\langle \Pis R_\varepsilon(E),E \rangle E+K, 
   \]
   where
   \[
    \|K(t)\|_F=\io(\|Y(t)-Y_\star\|_F^2)=\io(\|E(t)-E_\star\|_F^2).
   \]
   We recall that the expression of the orthogonal projection $\widehat{\Pi}_E$ onto the tangent space $\mathcal{T}_E\mathcal{S}_1$ at $E$ is given by
   \[
    \widehat{Z}:=\widehat{\Pi}_EZ=\Pis Z-\langle \Pis Z,E \rangle E, \qquad \forall Z\in \R^{n\times n}.
   \]
   In particular $\widehat{R}(E_\star)=0$ and
   \[
    \widehat{R}(E)=\widehat{R}(E)-\widehat{R}(E_\star)=\widehat{\Pi}_{E_\star} H_\varepsilon(E_\star) \widehat{\Pi}_{E_\star}(E-E_\star)+\io(\|E-E_\star\|_F^2),
   \]
   where $H_\varepsilon(E_\star)$ denotes the Hessian matrix of $F_\varepsilon$ at $E_\star$.
   Since it is assumed that $H_\varepsilon(E_\star)$ is positive definite and since 
   \[
    E-E_\star=\widehat{\Pi}_{E_\star}(E-E_\star)+\io(\|E-E_\star\|_F^2)
   \]
   we have, provided that $E$ is sufficiently close to $E_\star$,
   \[
    \frac{1}{2}\frac{d}{dt}\|E(t)-E_\star\|_F^2=\langle E-E_\star,\dot{E} \rangle=\langle E-E_\star,-\widehat{R}(E)+K \rangle=
   \]
   \[
    =\langle \widehat{\Pi}_{E_\star}(E-E_\star)+\io(\|E-E_\star\|_F^2),-\widehat{\Pi}_{E_\star} H_\varepsilon(E_\star) \widehat{\Pi}_{E_\star}(E-E_\star)+\io(\|E-E_\star\|_F^2) \rangle=
   \]
   \[
    =-\langle \widehat{\Pi}_{E_\star}(E-E_\star), H_\varepsilon(E_\star) \widehat{\Pi}_{E_\star}(E-E_\star) \rangle+\io(\|E-E_\star\|_F^3)\leq  
   \]
   \[
    \leq -\alpha \|\widehat{\Pi}_{E_\star}(E-E_\star)\|_F^2+\io(\|E-E_\star\|_F^3)\leq -\frac{\alpha}{2} \|E-E_\star\|_F^2,
   \]
   where $\alpha>0$ is the constant associated to the  strong minimum $E_\star$, that is
   \[
    \langle Z,H_\varepsilon(E_\star) Z\rangle \geq \alpha \|Z\|_F^2, \qquad \forall Z\in \mathcal{T}_{E_\star} \mathcal{S}_1.
   \]
   Hence $\|E(t)-E_\star\|_F$ decreases monotonically and exponentially to 0 as $t\rightarrow +\infty$.
   Similarly, since  $\dot{E}\in \mathcal{T}_E \mathcal{S}_1$,
   \[
    \frac{1}{\varepsilon}\frac{d}{dt}F_\varepsilon(E(t))=\langle R_\varepsilon(E),\dot{E}\rangle=\langle \widehat{R}_\varepsilon(E),\dot{E}\rangle=\langle \widehat{R}_\varepsilon(E),-\widehat{R}_\varepsilon(E)+K\rangle=
   \]
   \[
    =-\|\widehat{\Pi}_{E_\star}(E-E_\star), H_\varepsilon(E_\star) \widehat{\Pi}_{E_\star}(E-E_\star)\|_F^2+\io(\|E-E_\star\|_F^3)\leq 
   \]
   \[
    \leq-\alpha^2\|\widehat{\Pi}_{E_\star}(E-E_\star)\|_F^2+\io(\|E-E_\star\|_F^3)\leq -\frac{\alpha^2}{2}\|E-E_\star\|_F^2,
   \]
   which means that $F_\varepsilon(E(t))$ decreases monotonically and exponentially to $F_\varepsilon(E_\star)$ as $t\rightarrow +\infty$.

  \end{proof}

 \end{thm}
 
 Theorem \ref{thm_loc_conv} proves that, if a proper starting point is chosen, then the integration of the low rank ODE approaches a stationary point, as it happens for the gradient system \eqref{odeE}. Hence equation \eqref{odeY} can replace the original ODE. This fact can lead to computational benefit from the low rank underlying structure of the problem.

 \subsection{Implementation of the low rank \textit{inner iteration}}
 
 In this section we illustrate some details of the implementation of the \textit{inner iteration} for solving equation \eqref{prob_dk_nopen}, through a numerical integration of a system of ODEs. We show how to highlight the low rank properties of equation \eqref{odeY} by means of an equivalent system of ODEs and we discuss the numerical integration of the system. 
 
 Given an SVSD decomposition $Y(t)=USU^T$, we can rewrite equation \eqref{odeY} as
 \[
  \dot{U}S U^T+U\dot{S} U^T+US \dot{U}^T=-R+(I-UU^T)R(I-UU^T)+\eta US U^T. 
 \]
 Assuming $U^T\dot{U}=0$ yields the system
 \begin{equation}
  \label{odeUS}
  \begin{dcases}
   \dot{U}=-(I-UU^T)RU S^{-1}\\
  \dot{S}=-U^TRU+\eta S
  \end{dcases},
 \end{equation}
 which is equivalent to \eqref{odeY}. The matrix $S$ may lose the diagonal structure along the trajectory, but the SVSD decomposition $Y=USU^T$ still holds. System \eqref{odeUS} consists of two matrix ODEs of dimension $n$-by-$4$ and $4$-by-$4$ respectively. 
 
 Integration of  system \eqref{odeUS}  can be done in many ways. The simplest choice is the normalized explicit Euler method, which generally performs well. However in some cases the matrix $S$ may be close to singularity and Euler method may suffer the presence of the inverse of $S$ in its formulation. This problem can be overcome by means of a different integrator. Since we are not interested in the whole trajectory $Y(t)$, but only in the approximation of its stationary points, we can use a splitting method similar to that proposed in \cite{ceruti2022unconventional}. Algorithm \ref{alg_splitting} shows the outline of a single step integration of this approach.
 
 \begin{algorithm}
  \caption{Splitting method for the numerical integration step from $t_0$ to $t_1=t_0+h$}
  \label{alg_splitting}
  \begin{description}
   \item[Input:]  $U_0\in \R^{n\times 4}$ orthogonal and $S_0\in \textnormal{Sym}(\R^{4\times 4})$ non-singular such that $Y_0=U_0S_0U_0^T$
   \item[Output:] $U_1\in \R^{n\times 4}$ orthogonal and $S_1\in \textnormal{Sym}(\R^{4\times 4})$ invertible such that $Y_1=U_1S_1U_1^T$
  \end{description}
  \begin{algorithmic}[1]
   \State Begin
   \State Compute $K_1=U_0S_0+h\left(-R(Y_0)U_0+\langle P_{Y_0} R(Y_0),\Pis Y_0\rangle U_0 S_0\right)$
   \State Perform a QR factorization $K_1=U_1 T_1$ and compute $M_1=U_1^TU_0$
   \State  Define $\hat{S}_0=MS_0M^T$ and $\hat{Y}_0=U_1 \hat{S}_0 U_1^T$ \item Normalize $\hat{Y}_0$ and get $\tilde{Y}_0=U_1 \tilde{S}_0 U_1$ such that $\|\Pis\tilde{Y}_0\|_F=1$
   \State Compute $\tilde{R}_0=R(\tilde{Y}_0)$ and $\eta=\langle P_{\tilde{Y}_0}\tilde{R}_0,\Pis \tilde{Y}_0\rangle $
   \State Compute $\tilde{S}_1=\tilde{S}_0+hU_1^T\left(-\tilde{R}_0+\eta \tilde{Y}_0 \right)U_1$
   \State Normalize $\tilde{S}_1$ and get $S_1$ such that $\|\Pis (U_1 S_1 U_1^T)\|_F=1$
   \State Return $U_1$ and $S_1$
  \end{algorithmic}
 \end{algorithm}
 
 The choice of the stepsize is performed by means of an Armijo-type line search strategy as in \cite{guglielmi2022rank}, since the time derivative of the objective function is available. Provided a suitable starting point, theorem \ref{thm_loc_conv} ensures the convergence towards a stationary point. A possible choice for $U_0$ and $S_0$ comes from an SVSD decomposition of the gradient: this choice generally leads to a suitable approximation of a minimizer.
 We compute $x_0=x(L(W))$, $y_0=y(L(W))$ and $z_0=x_0\bullet x_0-y_0\bullet y_0$ and, since Remark \ref{rem_R} suggests to choose $Y_0=U_0S_0U_0^T=-R_0$, we define
 \begin{equation}
  \label{startUS}
  [U_0,D_0]=\texttt{qr}\left(\left(\begin{matrix}
   z_0+\one & z_0-\one & x_0 & y_0
  \end{matrix}\right),0\right), \qquad 
  S_0=-D_0\left(\begin{matrix}
   \frac{1}{4} & 0 & 0 & 0 \\
   0 & -\frac{1}{4} & 0 & 0 \\
   0 & 0 & -1 & 0 \\
   0 & 0 & 0 & 1 \\
  \end{matrix}\right)D_0^T,
 \end{equation}
 where $\texttt{qr}(\cdot,0)$ is the well-known \texttt{Matlab} function for the QR factorization. 
 However, during the \textit{outer iteration}, it could be more convenient to choose as starting value for the iteration of $\varepsilon_{l+1}$ the stationary points found in the $l$-th outer iteration, as  shown in Algorithm \ref{alg_outit}.
 
 \begin{algorithm}
  \caption{Inner iteration}
  \label{alg_innit}
  \begin{description}
   \item[Input:] A weight matrix $W$, a perturbation size $\varepsilon>0$, the starting values $U_0$ and $S_0$, an initial stepsize $h_0$, a tolerance tol and a maximum number of iterations maxit
   \item[Output:] The matrices $U_\star(\varepsilon)$ and $S_\star(\varepsilon)$ that form the solution of the optimization problem \eqref{prob_ii} $E_\star(\varepsilon)=\Pis(U_\star(\varepsilon)S_\star(\varepsilon)U_\star(\varepsilon)^T)$
  \end{description}
  \begin{algorithmic}[1]
   \State Begin
   \State Initialize $U_0$ and $S_0$ (e.g. by means of \eqref{startUS})
   \State Compute $f_0=F(\Pis (U_0 S_0 U_0^T))$ and set $f_1=f_0+1$ 
   \State Set $j=0$ 
   \While{$|f_1-f_0|>\textnormal{tol}$ \textnormal{and} $j<\textnormal{maxit}$}
   \State With an Armijo stepsize choice, perform Algorithm \ref{alg_splitting} and compute $U_1$ and $S_1$.   
   \State Update $f_0=f_1$,  $f_1:=F(\Pis(U_1 S_1 U_1^T))$ and set $j:=j+1$
   \EndWhile
  \end{algorithmic}
 \end{algorithm}

 \section{The \textit{outer iteration}}
 
 Once that a computation of the optimizers is available for a given $\varepsilon>0$ and a fixed $k$, we need to determine an optimal value for the perturbation size.
 Let $E_\star(\varepsilon)$ be a solution of the optimization problem \eqref{prob_ii} and consider the function
 \[
  \varphi(\varepsilon)=F_\varepsilon(E_\star(\varepsilon)).
 \]
 This function is non-negative and we define $\varepsilon_\star$ as the smallest zero of $\varphi$. Assuming that the $k$-th and $(k+1)$-st eigenvalues of $L(W+\varepsilon E_\star(\varepsilon))$ are simple, for $0\leq \varepsilon < \varepsilon_\star$, yields that $\varphi$ is a smooth function in the interval $[0,\varepsilon_\star)$.
 The aim of the \textit{outer iteration} is to approximate $\varepsilon_\star$, which is the solution of the optimization problem \eqref{prob_dk_nopen}.
 
 In order to find $\varepsilon_\star$ we use a combination of the well-known Newton and bisection methods, which provides an approach similar to \cite{guglielmi2015low,guglielmi2016method} or \cite{guglielmi2022rank}. If the current approximation $\varepsilon$ lies in a left neighborhood of $\varepsilon_\star$, it is possible to exploit Newton's method, since $\varphi$ is smooth there; otherwise in a right neighborhood we use the bisection method (see Algorithm \ref{alg_outit}). The following result provides a simple formula for the first derivative of $\varphi$ required by Newton's method.
 
 \begin{lem}
  \label{lem_derfeps}
  It holds
  \[
   \varphi'(\varepsilon)=\frac{d}{d\varepsilon} F_\varepsilon(E_\star(\varepsilon))=\langle G_\varepsilon(E_\star(\varepsilon)),E_\star(\varepsilon)\rangle=-\|G_\varepsilon(E_\star)\|_F.
  \]
  \begin{proof}
   As shown for the time derivative formula, we get
   \[
    \varphi'(\varepsilon)=\frac{d}{d\varepsilon} \left(\lambda_{k+1}(\varepsilon L(E_\star(\varepsilon)))- \lambda_k( \varepsilon L(E_\star(\varepsilon)))\right)=
   \]
   \[
    =x^T\frac{d}{d\varepsilon} (\varepsilon L(E_\star(\varepsilon)))x-y^T\frac{d}{d\varepsilon} (\varepsilon L(E_\star(\varepsilon)))y=\langle xx^T-yy^T,L(E_\star(\varepsilon))+\varepsilon(L(E'_\star(\varepsilon)))\rangle=
   \]
   \[
    =\langle L^*(x(\varepsilon)x(\varepsilon)^T-y(\varepsilon)y(\varepsilon)^T),E_\star(\varepsilon)\rangle=\langle G_\varepsilon(E_\star(\varepsilon)),E_\star(\varepsilon)+\varepsilon E'_\star(\varepsilon)\rangle.
   \]
   where $E'_\star(\varepsilon)$ is the derivative with respect to $\varepsilon$ of $E_\star(\varepsilon)$. 
   Since $E_\star$ is a unit norm stationary point of \eqref{prob_ii} and a zero of the derivative of the objective functional $F_\varepsilon$, then $G_\varepsilon(E_\star(\varepsilon))$ is a negative multiple of $E_\star$. Thus 
   \[
    G_\varepsilon(E_\star)=-\|G_\varepsilon(E_\star)\|_F E_\star(\varepsilon), \qquad
    \langle G_\varepsilon(E_\star(\varepsilon)), E'_\star(\varepsilon)\rangle=\frac{\langle E_\star(\varepsilon), E'_\star(\varepsilon) \rangle}{\langle G_\varepsilon(E_\star(\varepsilon)), E_\star(\varepsilon)\rangle}=0.
   \]

  \end{proof}
 \end{lem}
 
 \begin{algorithm}
  \caption{Outer iteration}
  \label{alg_outit}
  \begin{description}
   \item[Input:] A weight matrix $W$, an interval and an initial guess $\varepsilon_0\in[\varepsilon_{\textnormal{lb}},\varepsilon_{\textnormal{ub}}]$ for $\varepsilon_\star$, a tolerance toler and a maximum number of iterations niter
   \item[Output:] The value $\varepsilon_\star$ and the minimizer associated $E_\star(\varepsilon_\star)$
  \end{description}
  \begin{algorithmic}[1]
   \State Begin
   \State Perform Algorithm \ref{alg_innit} and compute $E_\star(\varepsilon_0)$  
   \State Set $l=0$
   \While{$l<\textnormal{niter}$ \textnormal{and} $\varepsilon_{\textnormal{ub}}-\varepsilon_{\textnormal{lb}}>\textnormal{toler}$}
   \If{$\varphi(\varepsilon_l)<\textnormal{toler}$}
   \State Set $\varepsilon_{\textnormal{ub}}:=\min(\varepsilon_{\textnormal{ub}},\varepsilon_l)$
   \State Set $\varepsilon_{l+1}:=\frac{\varepsilon_{\textnormal{lb}}+\varepsilon_{\textnormal{ub}}}{2}$ (bisection step)
   \Else
   \State Set $\varepsilon_{\textnormal{lb}}:=\max(\varepsilon_{\textnormal{lb}},\varepsilon_l)$
   \State Compute $\varphi(\varepsilon_l)$ and $\varphi'(\varepsilon_l)$ 
   \State Update $\varepsilon_{l+1}:=\varepsilon_l-\frac{\varphi(\varepsilon_l)}{\varphi'(\varepsilon_l)}$ (Newton step)
   \EndIf
   \If{$\varepsilon_{l+1}\notin[\varepsilon_{\textnormal{lb}},\varepsilon_{\textnormal{ub}}]$}
   \State Set $\varepsilon_{l+1}:=\frac{\varepsilon_{\textnormal{lb}}+\varepsilon_{\textnormal{ub}}}{2}$
   \EndIf   
   \State Set $l:=l+1$
   \State Compute $E_\star(\varepsilon_l)$ by applying Algorithm \ref{alg_innit} with starting value $E_\star(\varepsilon_{l-1})$
   \EndWhile
   \State Return $\varepsilon_\star:=\varepsilon_l$ and $E_\star(\varepsilon_\star)$
  \end{algorithmic}
 \end{algorithm}
  
 Finally we perform Algorithm \ref{alg_outit} for some values of $k\in [k_{\min},k_{\max}]$ and we select the index of the largest structured distance computed.
 
 \subsection{The penalized version}
 
 A similar approach can be followed for the penalized problem \eqref{prob_ii_pen}. For $\varepsilon,c>0$, let $E_\star(\varepsilon,c)$ be a solution of the penalized \textit{inner iteration} \eqref{prob_ii_pen} and consider the function
 \[
  \varphi_c(\varepsilon)=F_{\varepsilon,c}(E_\star(\varepsilon,c))
 \]
 and define $\varepsilon_\star$ as the minimum zero of $\varphi_c$. Again assuming that the $k$-th and $(k+1)$-st eigenvalues of $L(W+\varepsilon E(\varepsilon,c))$ are simple, for $0\leq \varepsilon < \varepsilon_\star$, yields that $\varphi_c$ is a smooth function in the interval $[0,\varepsilon_\star)$.
 \begin{lem}
  \label{lem_derfeps_pen}
  It holds
  \[
   \varphi_c'(\varepsilon)=\frac{d}{d\varepsilon} F_{\varepsilon,c}(E_\star(\varepsilon))=\langle G_{\varepsilon,c}(E_\star(\varepsilon,c)),E_\star(\varepsilon,c)\rangle=-\|G_{\varepsilon,c}(E_\star(\varepsilon,c))\|_F.
  \]
  \begin{proof}
   It is a direct consequence of Lemma \ref{lem_derfeps}.
  \end{proof}
 \end{lem}
  
 As done for the low rank method, it is possible to implement an algorithm to solve problem \eqref{prob_dk} by introducing the penalization term. The only difference is in the \textit{inner iteration} where we integrate equation \eqref{odeE_pen}, for instance with normalized Euler method.

 \section{Numerical experiments}
 
 In this section we compare the behavior of the spectral gaps 
 \[
  g_k(W)=\lambda_{k+1}(W)-\lambda_k(W)
 \]
 and the structured distance to ambiguity as stability indicators. For the computation of $d_k(W)$ we use both Algorithm \ref{alg_outit} and the integration of the full-rank system. To distinguish between these two results, we will denote, respectively, the unstructured  distances as
 \[
  d_k^{LOW}(W), \qquad d_k^{FULL}(W).
 \]
 If the optimizer found is not admissible, we integrate only the penalized equation \eqref{odeE_pen}, since the low rank system is not suitable.
 
 We present four different examples with different features: in the first three the penalization term is not required and hence we can use Algorithm \ref{alg_outit}, while the last one shows a non-common case where the non-negativity constraint must be taken into account. In all experiments we set the tolerance of $10^{-2}$ for the \textit{outer iteration}.
 
 \subsection{A slightly sparse example: the JOURNALS matrix}
 
\begin{figure}
  \begin{center}
   \includegraphics[scale=0.5]{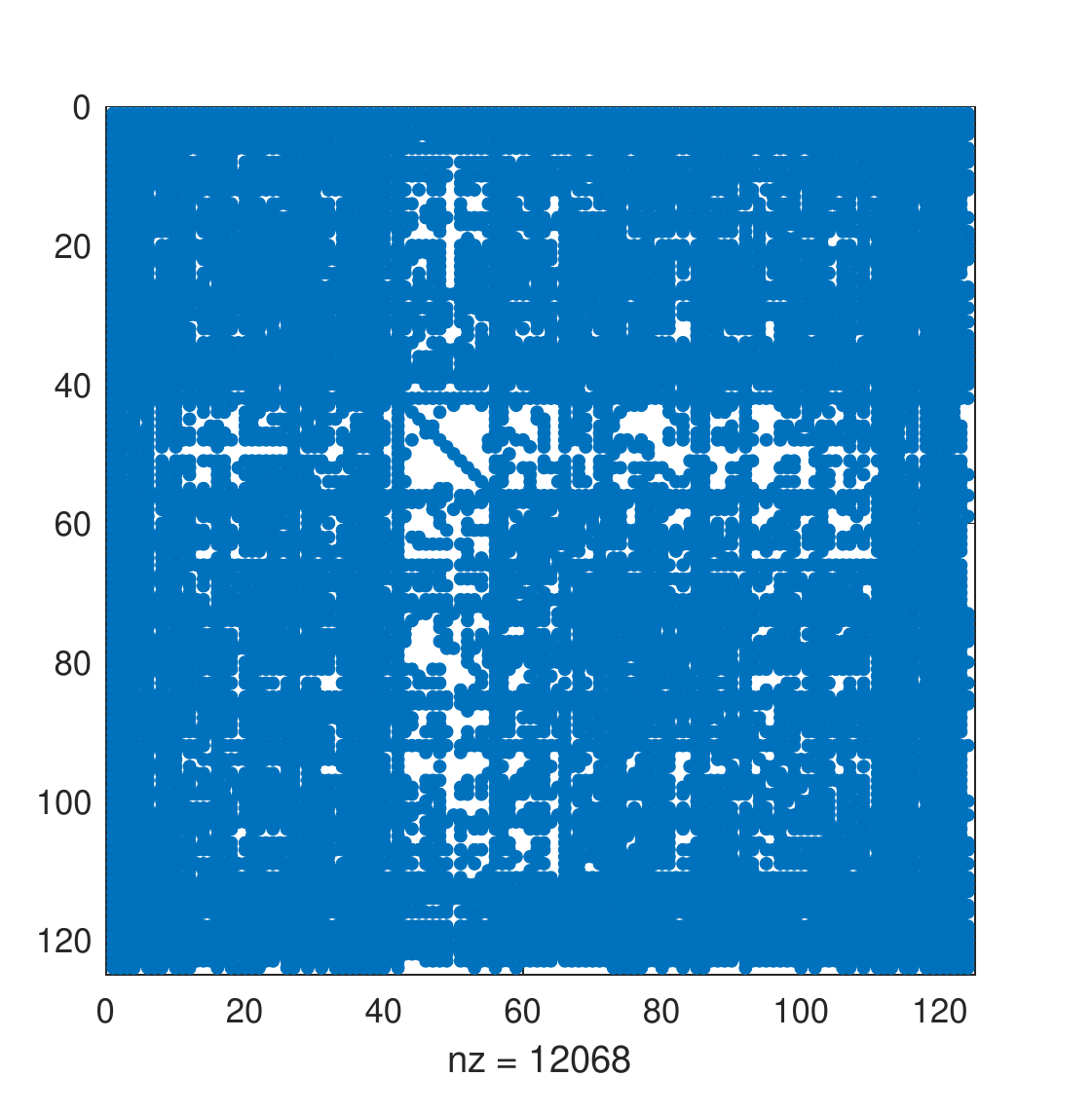}
   \caption{Structural pattern of the JOURNALS matrix.}
  \end{center}
  \label{fig_Journals}
 \end{figure}

 The JOURNALS matrix comes from a Pajek network converted to a sparse adjacency matrix for inclusion in the University of Florida SuiteSparse matrix collection (see   \cite{davis2011university} for more details). It represents an undirected weighted graph with $n=124$ vertices and $m=12068\approx 97n$ edges, whose structural pattern is shown in figure \ref{fig_Journals}. The pattern of $W$ does not suggest a suitable number of clusters to partition the graph. We select $k\in[3,8]$ and we compare in Table \ref{tab_Journals} the results of the unstructured distances with the spectral gaps, where we have set the inner tolerance as $10^{-4}$. In this case the two criteria for the choice of the best number of cluster disagree: while the structured distances select $k_\textnormal{opt}=8$, the unstructured distance, i.e. the largest spectral gap, prefers $k=4$. In this example, the size and the pattern of the matrix implies that the low rank ODE is more convenient than the full rank gradient system. In particular the gain in memory requirement is given by the ratio between $m$ (for $E$ of ODE \eqref{odeE}) and $4n+16$ (for $U$ ans $S$ of system \eqref{odeUS}), that is
 \[
  \frac{m}{4n+16}=\frac{12068}{4 \cdot 124+16}\approx 23.57,
 \]
 which is quite convenient. However the CPU time of the low rank system method is $49$ seconds against the $24$ seconds of the full rank system. A possible reason behind this behavior is that the gradient system requires less effort to reach convergence and hence less eigenvalues computations, which are the most expensive procedures in all the methods.

 \begin{table}
  \centering
  \[  
  \begin{array}{|c|c|c|c|c|}
   \hline
   k & g_k(W) & d_k^{LOW}(W)  &  d_k^{FULL}(W) & \left|d_k^{LOW}(W)-d_k^{FULL}(W)\right| \\
   \hline
   3 & 8.9963 & 5.3128 &  5.3128  & <10^{-4} \\
   \hline
   4  &  \textbf{39.1923} &  6.3225 &     6.3195 & 0.0030 \\
   \hline
   5 & 27.9853 & 5.3637 & 5.3638 & 0.0001 \\
   \hline
   6 & 10.4987 & 4.4254 &    4.5364 &  0.1109 \\
   \hline
   7 & 22.0178 & 7.4290 &   7.4296 & 0.0006 \\
   \hline
   8 & 33.9873 &  \textbf{8.1434} &    \textbf{8.1652} &  0.0218 \\
   \hline
  \end{array}
  \]
  \caption{Comparison with the distances for JOURNALS matrix. The highlighted results indicate the best value for $k$ according to each method.}
  \label{tab_Journals}
 \end{table}
 
\subsection{A Machine Learning example: the ECOLI matrix}

 The ECOLI matrix is a Machine Learning dataset from the SuiteSparse Matrix Collection and the UCI Machine Learning Repository (see \cite{pasadakis2022multiway}) that describes the protein localization sites of the bacteria \textit{E. coli}.   
 For $n$ data points, the connectivity matrix $C\in \R^{n\times
 n}$ is created from a $k$-nearest neighbors routine, with $k$ set such that
 the resulting graph is connected. The similarity matrix $S \in         
 \R^{n\times n}=(s_{ij})$ between the data points is defined as          
 \[                                                      
  s_{ij} = \max\{s_i(j), s_j(i)\} \;\; \text{with}\;                 
  s_i(j) = \exp \left(-4 \frac{\|x_i - x_j \|^2}{\sigma_i^2} \right)            
 \]                                                                                                                   
 with $\sigma_i$ standing for the Euclidean distance between the $i$-th  
 data point and its closest $k$-nearest neighbor. The adjacency matrix $W$
 is then created as $W=C\bullet S$ (here $\bullet$ denotes the componentwise product). 
 
\begin{figure}
  \begin{center}
   \includegraphics[scale=0.45]{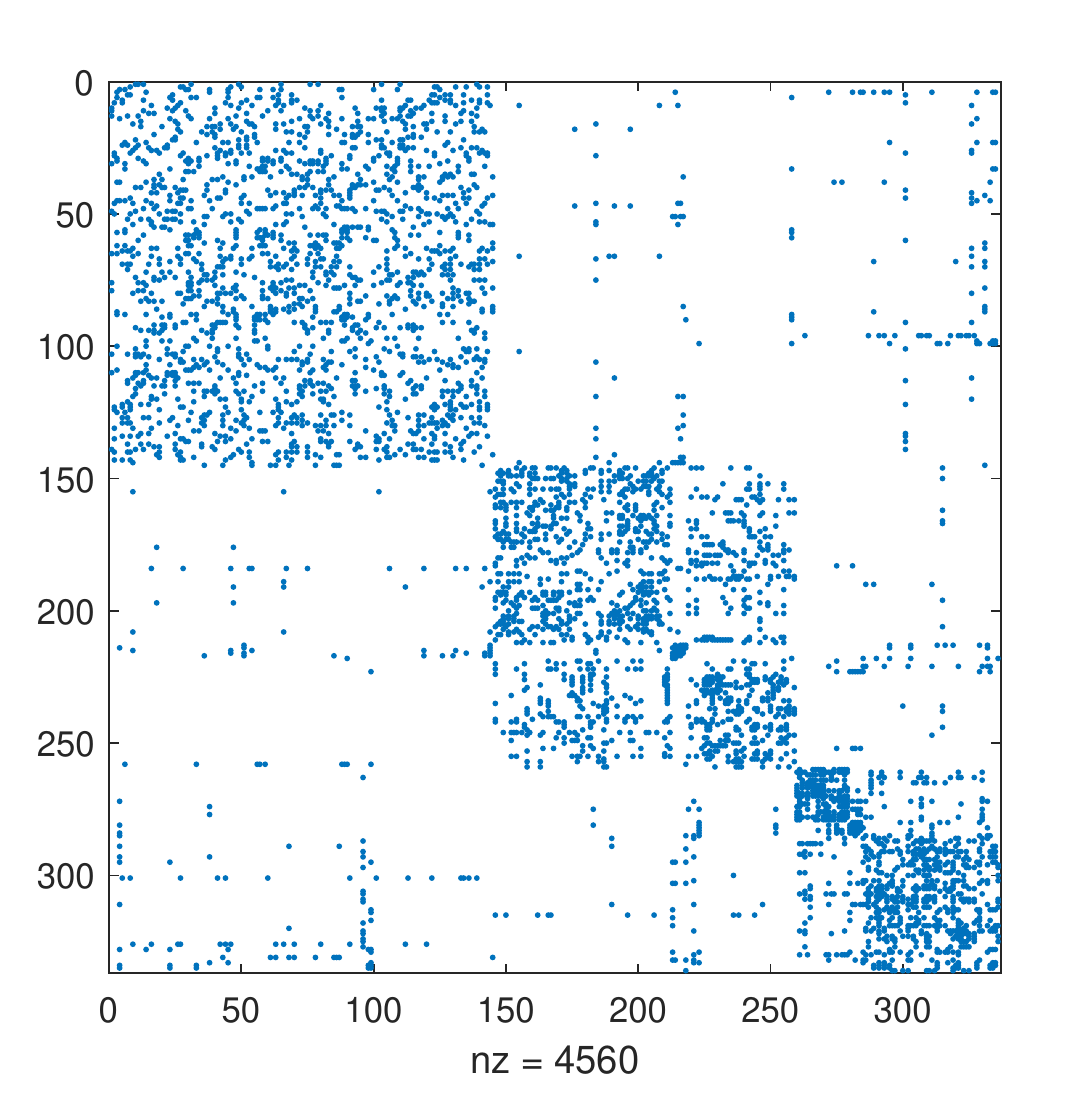}
   \hspace{1.5cm}
   \includegraphics[scale=0.45]{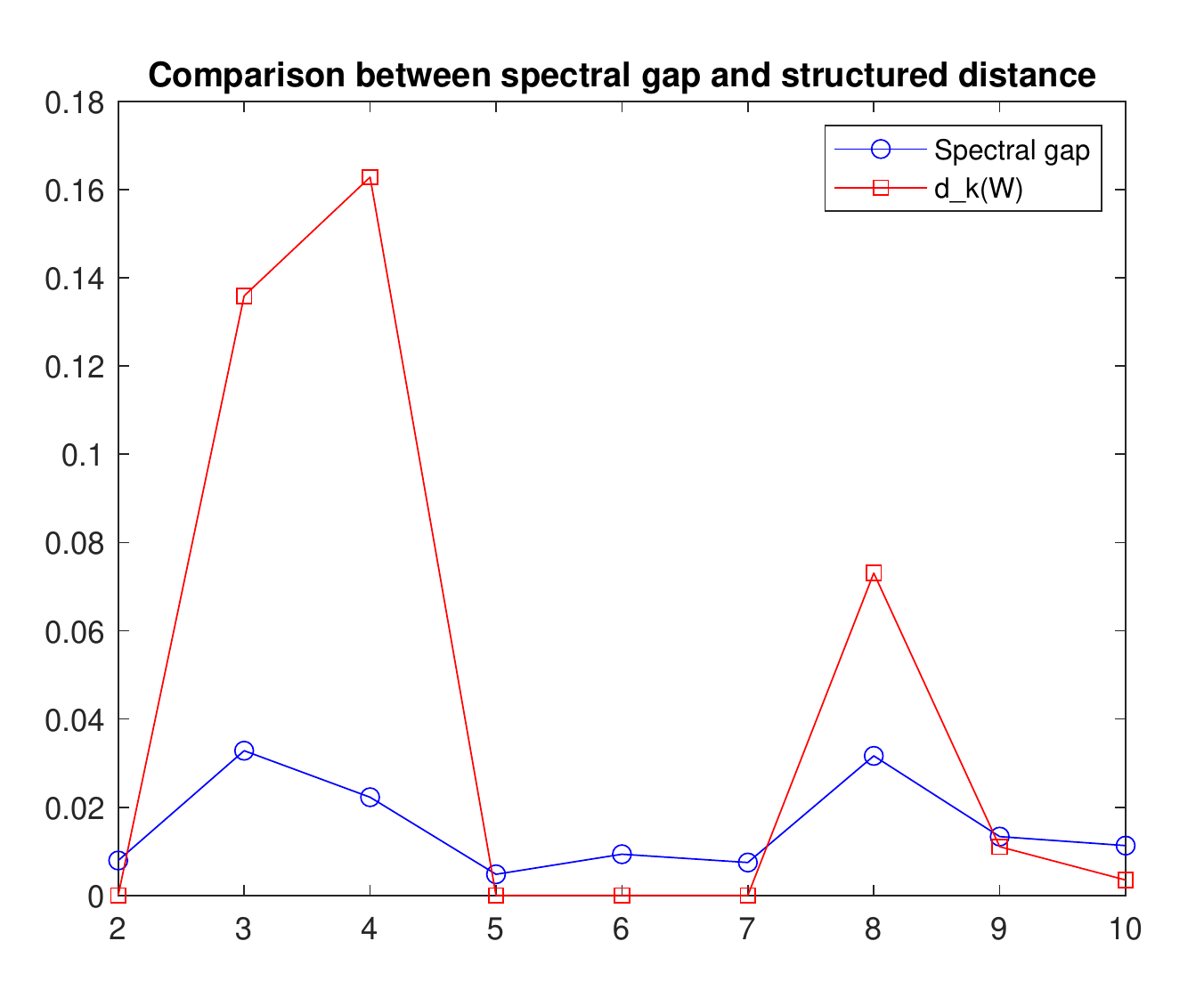}
   \caption{ECOLI matrix: on the left its structural pattern, on the right the spectral gap and $d_k(W)$}
  \end{center}
  \label{fig_Ecoli}
 \end{figure}

 This matrix has $n=336$ vertices and $m=4560\approx 13.6n$ edges and its pattern is shown in figure \ref{fig_Ecoli}. In this case the structure of $W$ contains three possible clusters and the second and third appear to be split in two sub-communities. This facts suggest that a suitable number of clusters should be $k\in\{3,4,5\}$. We will test the capability of all the methods to identify this feature with inner tolerance $10^{-9}$. The graph in Figure \ref{fig_Ecoli} shows the performances of the methods. In this case the low rank and full rank unstructured distances coincide, up to machine errors, and thus they are reported as $d_k(W)$. It is evident that for all the methods the best choices are $k=3,4$ or $8$. More precisely the largest spectral gap is $g_3(W)=0.328$ slightly greater than $g_8(W)=0.316$, while for the unstructured distance the best choice is $k=4$, a bit preferable than $k=3$. Also in this case there is a gain in memory saving 
 \[
  \frac{m}{4n+16}=\frac{4560}{4 \cdot 336+16}\approx 3.35,
 \]
 while the ratio in CPU time is almost the same of the previous example: $18$ seconds for the low rank while $8$ seconds for the gradient system.
 
 \subsection{A social network community: the EGO-FACEBOOK matrix}
 
\begin{figure}
  \begin{center}
   \includegraphics[scale=0.4]{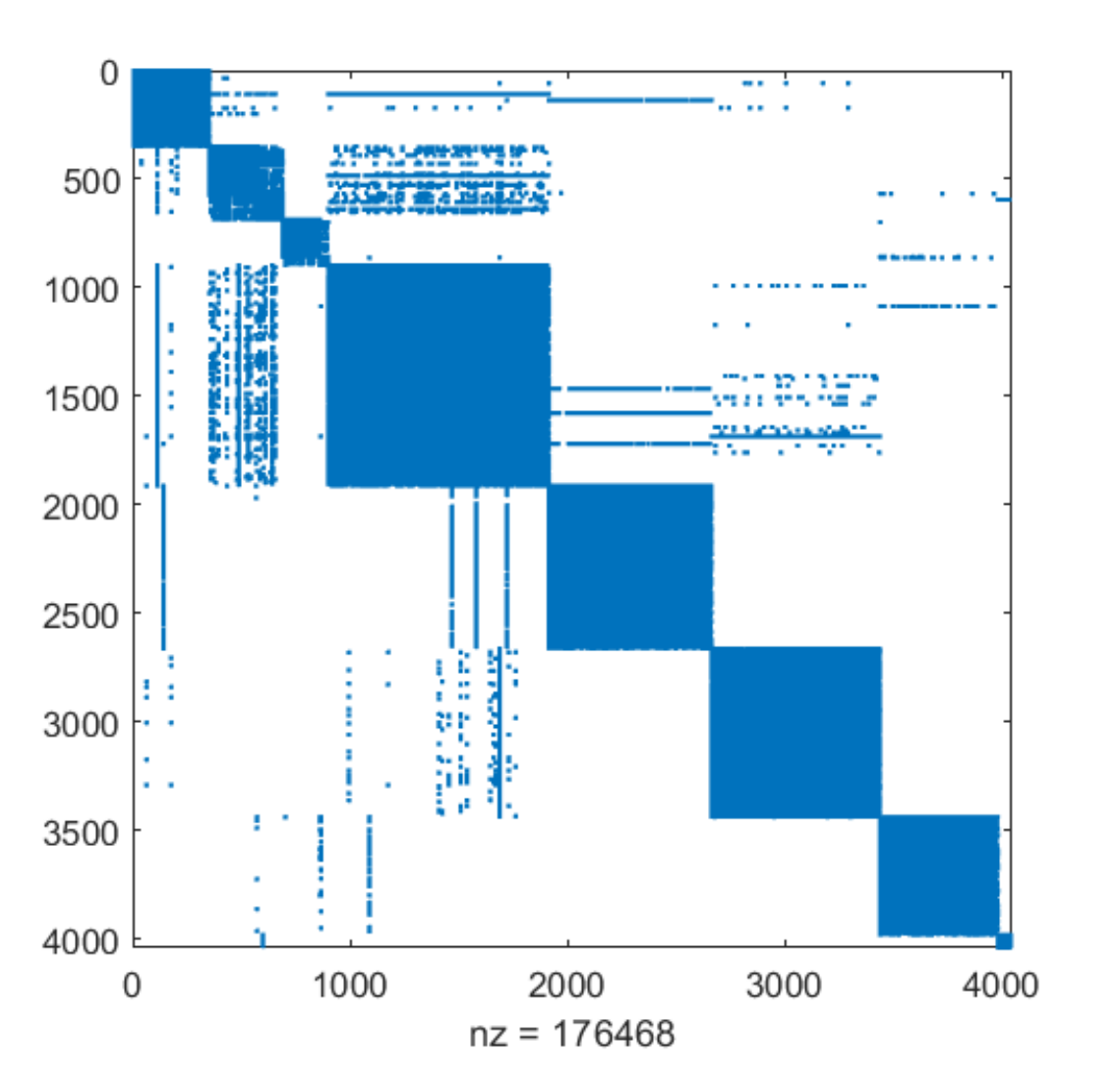}
   \includegraphics[scale=0.4]{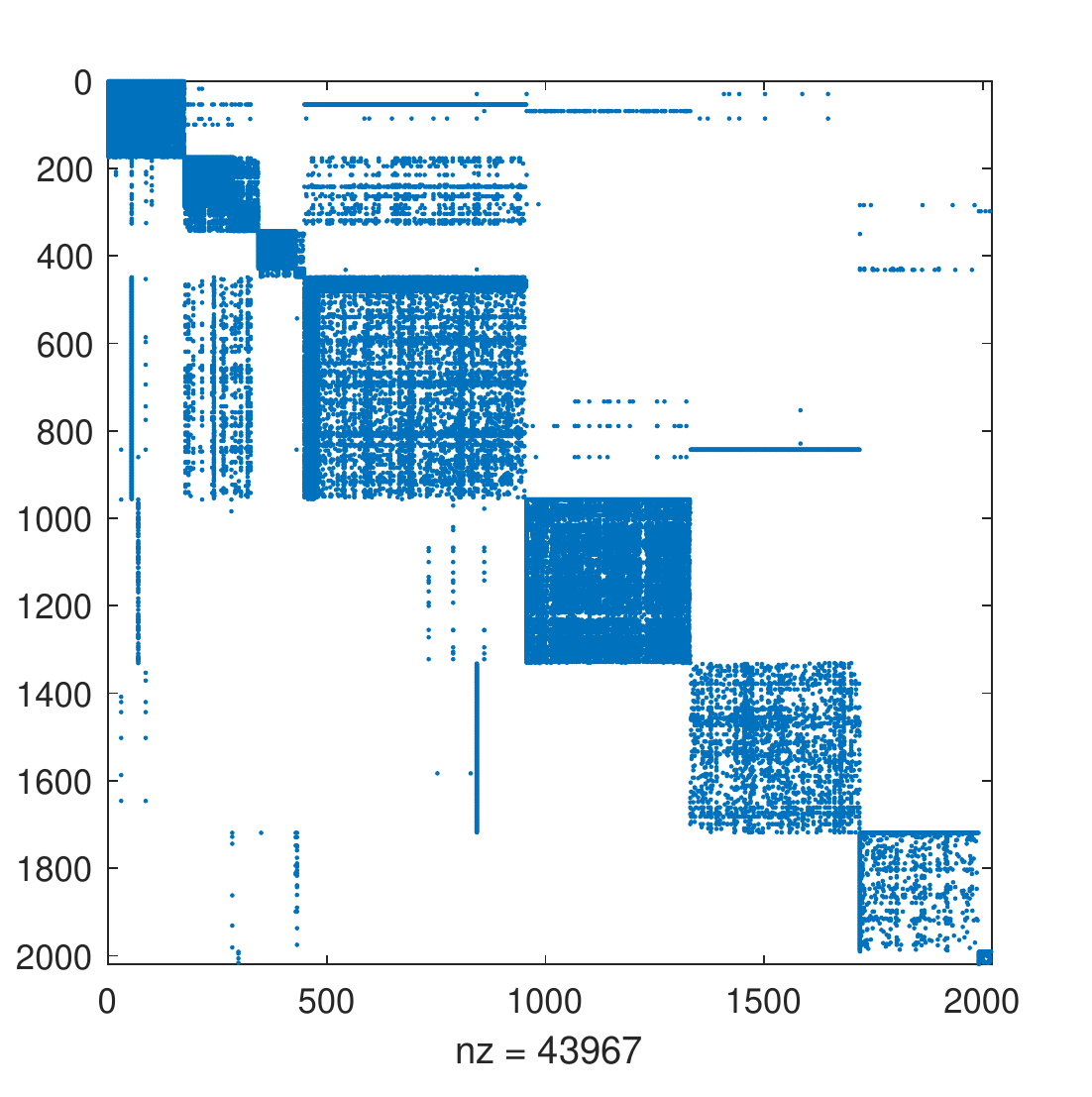}
   \includegraphics[scale=0.4]{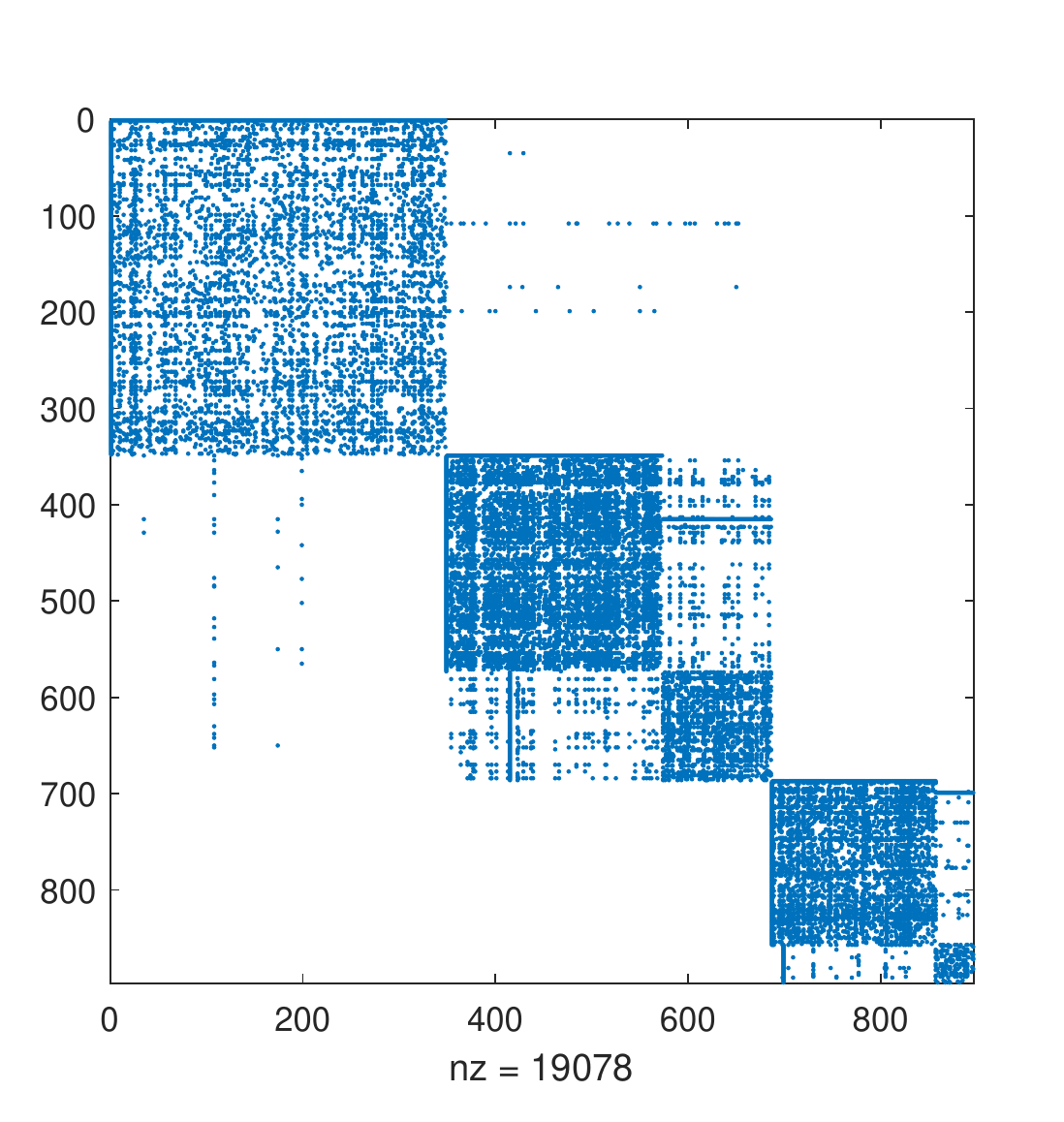}
   \caption{EGO-FACEBOOK matrices: on the left the whole structural pattern of the full matrix $W_1$, in the middle the compressed matrix $W_2$ and on the right the sub-matrix $W_3$.}
  \end{center}
  \label{fig_ego-Fb_patterns}
 \end{figure}

 The EGO-FACEBOOK matrix represents a dataset that consists of ``circles'' (or ``friends lists'') of the social network Facebook from the SNAP dataset (see \cite{leskovec2012learning} for more details). The data were collected from survey participants using this Facebook app. The whole matrix $W_1$ has $n_1=4093$ vertices with $m_1=176468\approx 43.1 n_1$ edges (see Figure \ref{fig_ego-Fb_patterns}). In order to perform further tests of the algorithms, we also consider two reduced versions, $W_2$ and $W_3$, of the whole matrix $W_1$. 
 
 The matrix $W_2$ is obtained by means of a compression that mantains the pattern and the density of the original matrix, but it halves the dimension: more precisely we define
 \[
  (W_2)_{i,j}=\frac{(W_1)_{2i-1,2j-1}+(W_1)_{2i-1,2j}+(W_1)_{2i,2j-1}+(W_1)_{2i,2j}}{4}, \qquad i,j=1,\dots, \frac{n_1-1}{2}
 \]
 and then we set to zero the entries with the smallest value such that the compressed matrix has the same density of $W_1$. We obtain the matrix $W_2$ with $n_2=2019$ vertices and $m_2=43967\approx 21.8 n_2$ edges. Finally we considered the main minor $W_3$ of $W_1$ formed by the first $n_3=896$ vertices and with $m_3=19078\approx 21.3 n_3$ edges, which contains the first three main blocks of the whole matrix.
 
 \subsubsection{Whole matrix}
 
 First we analyze the whole matrix $W_1$. In Table \ref{tab_Fb_full} we report the results where we set the inner tolerance $10^{-9}$, which will be also the accuracy of the experiments for $W_2$ and $W_3$. The criteria disagree also in this case: the unstructured distance prefers $k=6$, while the largest spectral gap is $k=7$. The factor for the memory saving gain with respect to the full rank system is
 \[
  \frac{m_1}{4n_1+16}=\frac{176468}{4 \cdot 4093+16}\approx 10.77,
 \]
 while the CPU time performances are $1270$ seconds for the low rank and $424$ seconds for the gradient system.
 
 \begin{table}
  \centering
  \[  
  \begin{array}{|c|c|c|c|c|}
   \hline
   k & g_k(W_1) & d_k^{LOW}(W_1) &  d_k^{FULL}(W_1) & \left|d_k^{LOW}(W_1)-d_k^{FULL}(W_1)\right| \\
   \hline
   3 & 0.0182 & 1.1015 &    1.1015 & <10^{-4} \\
   \hline
   4 & 0.0211 & 3.2270 &   3.2270 & <10^{-4} \\
   \hline
   5 & 0.0423 & 5.7524 &    5.7524 & 0 \\
   \hline   
   6 & 0.0526 & \textbf{6.6343} &    \textbf{6.6343} & 0 \\
   \hline 
   7 & \textbf{0.5153} & 1.1798 &   1.1798 & <10^{-4} \\
   \hline 
   8 & 0.0546 & 0.4725 &   0.4850 & 0.0125 \\
   \hline
  \end{array}
  \]
  \caption{Comparison between the distances for the full EGO-FACEBOOK matrix $W_1$. The marked bold results indicate the best value for $k$ according to each method.}
  \label{tab_Fb_full}
 \end{table}

 \subsubsection{Compressed matrix}

 In order to test the robustness of the algorithms, we compare the results between the full matrix $W_1$ and its compressed version $W_2$. Table \ref{tab_Fb_compressed} shows that the algorithms gives the same optimal values of the whole matrix computation, even though there are some differences in magnitude. The factor for the memory saving gain with respect to the full rank system is
 \[
  \frac{m_2}{4n_2+16}=\frac{43967}{4 \cdot 2019+16}\approx 5.43,
 \]
 while the CPU time performances are $147$ seconds for the low rank and $66$ seconds for the gradient system. This means that the computational time has scaled approximately of a factor $8$ between $W_1$ and $W_2$.
 
 \begin{table}
  \centering
  \[  
  \begin{array}{|c|c|c|c|c|}
   \hline
   k & g_k(W_2) & d_k^{LOW}(W_2) &  d_k^{FULL}(W_2) & \left|d_k^{LOW}(W_2)-d_k^{FULL}(W_2)\right| \\
   \hline
   3 & 0.0037 & <10^{-7} & <10^{-7} & 0 \\
   \hline
   4 & 0.0155 & 0.2996 &   0.2996 & <10^{-4} \\   \hline
   5 & 0.0030 & <10^{-7}& <10^{-7} & 0 \\
   \hline
   6 & 0.0638 & \textbf{1.7965} &    \textbf{1.7613} & 0.0352 \\
   \hline
   7 & \textbf{0.3865} & 1.2727 &    1.2727 & <10^{-4} \\
   \hline
   8 & 0.0036 & <10^{-7} & <10^{-7} & 0 \\
   \hline
  \end{array}
  \]
  \caption{Comparison between the distances for the compressed EGO-FACEBOOK matrix $W_2$. The marked bold results indicate the best value for $k$ according to each method.}
  \label{tab_Fb_compressed}
 \end{table}

 \subsubsection{Reduced matrix}
 
 Now we focus on $W_3$. From its pattern it is clear that the most reasonable choices for the number of clusters $k$ should be between $3,4$ and $5$. We also include $k=6$ and we investigate for that values  the performances of the methods. In this case all the methods agree and the best number of clusters is $k=3$. The factor for the memory saving gain with respect to the full rank system is
 \[
  \frac{m_3}{4n_3+16}=\frac{19078}{4 \cdot 896+16}\approx 5.30,
 \]
 while the CPU time performances are $34$ seconds for the low rank and $17$ seconds for the gradient system. 
 
 \begin{table}
  \centering
  \[  
  \begin{array}{|c|c|c|c|c|}
   \hline
   k & g_k(W_3) & d_k^{LOW}(W_3)  &  d_k^{FULL}(W_3) & \left|d_k^{LOW}(W_3)-d_k^{FULL}(W_3)\right| \\
   \hline
   3 &  \textbf{0.6040} &  \textbf{1.4092} &  \textbf{1.4795}  &  0.0703 \\
   \hline
   4  &  0.1724  &  0.7529  &  0.7529  & 0 \\
   \hline
   5 &  0.1870  &   0.2650 &    0.2629  &  0.0020\\
   \hline
   6 & <10^{-7} & <10^{-7} & <10^{-7} & 0\\
   \hline
  \end{array}
  \]
  \caption{Comparison between the distances for EGO-FACEBOOK matrix reduced $W_3$. The highlighted results indicate the best value for $k$ according to each method.}
  \label{tab_Fb_reduced}
 \end{table}

 \subsection{An example with penalization: the Stochastic Block Model (SBM)}
 
\begin{figure}
  \begin{center}
   \includegraphics[scale=0.5]{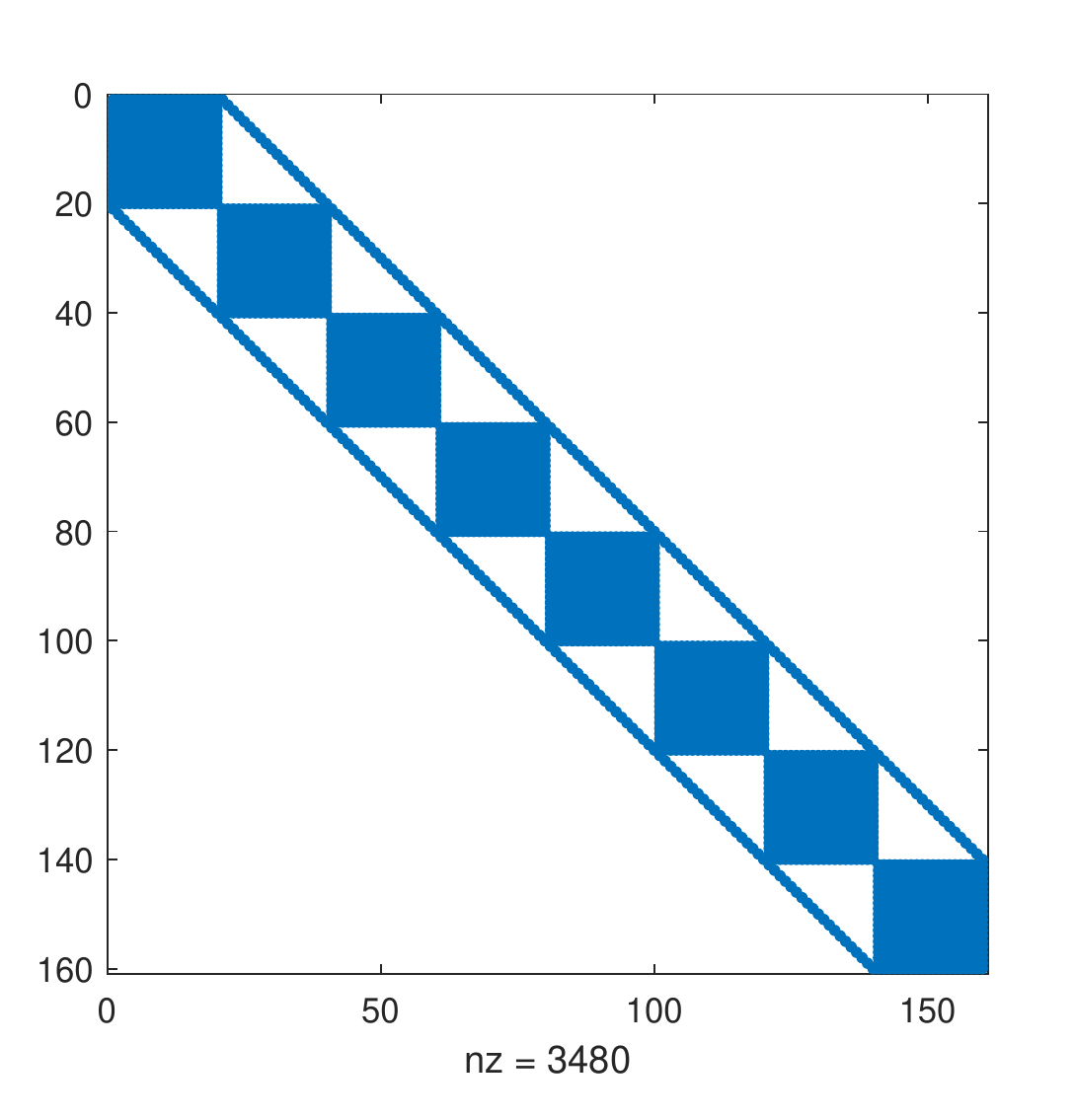}
   \caption{SBM matrix structural pattern.}
  \end{center}
  \label{fig_SBM_pattern}
 \end{figure}

 The Stochastic Block Model (SBM) is a model of generating random graphs that tend to have communities. It is an important model in a wide range of fields, from sociology to physics. In this example we consider $n=160$ vertices partitioned in $p=8$ clusters $C_1,\dots,C_p$ of $q=20$ elements each. We consider a random full symmetric matrix $J\in \R^{q\times q}$ and we build the matrix
 \[
  W=\textnormal{kron}(I_p,J)+\textnormal{kron}(B_p,I_q), \qquad
  B_p=\begin{pmatrix}
       0 & 1 & 0 & \dots & 0 \\
       1 & 0 & 1 & \ddots & \vdots \\
       0 & \ddots & \ddots & \ddots & 0 \\
       \vdots & \ddots & 1 & \ddots & 1\\
       0 & \dots & 0 & 1 & 0 \\
      \end{pmatrix}\in \R^{p\times p},
 \]
 where kron denotes the Kronecker product.
 The weight matrix generated has the pattern in Figure \ref{fig_SBM_pattern}, with $m=3480\approx 21.75 n$ non-zero entries and it has $p$ blocks by construction. If we apply Algorithm \ref{alg_outit}, some values of $k$ provide a non-admissible solution, which means that in this case a penalization is needed and the low rank system \eqref{odeUS} cannot be exploited. In particular for $k=8$, which is one of the candidate optimal values, the non-negativity constraint break cannot be ignored. In Table \ref{tab_SBM} we show the results of the integration of the full rank gradient system \eqref{odeE_pen}, where we introduce in the $j$-th \textit{inner iteration} a penalization $c_j$ that starts from $c_0=0$ and then increases by adding $0.5$ during each iteration, that is $c_j=0.5j$. The results found are admissible or slightly not, with the norm of the negativity part that is of order $10^{-5}$: in the last case we ensure that the optimizer is admissible by removing this error.
  The time required by the computation is $14$ seconds.
  
 \begin{table}
  \centering
  \[  
  \begin{array}{|c|c|c|c|c|}
   \hline
   k & g_k(W) &   d_k^{FULL}(W) \\
   \hline
   3 & 0.6488 & 6.2536 \\
   4 & 0.7654 & 8.1194 \\
   5 & 0.7654 & 6.2051 \\
   6 & 0.6488 & 8.3657 \\
   7 & 0.4335 & 5.7874 \\
   8 & \textbf{12.1650}  & \textbf{12.4264} \\
   9 & 0.1522 & 0.7894\\
   \hline
  \end{array}
  \]
  \caption{Comparison with the distances for SBM matrix. The marked bold results indicate the best value for $k$ according to each method.}
  \label{tab_SBM}
 \end{table}

 \nocite{*}
 \bibliographystyle{plain}
 \bibliography{graph_ode_article}

\end{document}